\documentclass[11pt]{article}

\usepackage{epsf,epsfig,amsmath,amsfonts,amssymb}
\usepackage{delarray}
\usepackage{mathrsfs}
\usepackage{enumerate}
\usepackage[colorlinks,
            linkcolor=purple,
            anchorcolor=brown,
            citecolor=blue
            ]{hyperref}
\usepackage{xcolor}

\newtheorem{theorem}{Theorem}
\newtheorem{example}{Example}
\newtheorem{corollary}[theorem]{Corollary}
\newtheorem{proposition}[theorem]{Proposition}
\newtheorem{lemma}[theorem]{Lemma}
\newtheorem{definition}{Definition}
\newtheorem{remark}{Remark}

\newcommand\goldenbox{\setlength\fboxrule{0.02em}\setlength\fboxsep{0em}%
\fbox{\phantom{\rule{1ex}{1.618ex}}}}

\def\endofproof{\hbox{}\hfill\goldenbox}

\begin{document}

\title{\sc Martingale Problem under Nonlinear Expectations}

\author{
{\sc Xin Guo\footnote{Department of Industrial Engineering and
Operations Research, UC Berkeley, CA 94720-1777, USA,
\texttt{xinguo@ieor.berkeley.edu} }\ \,
and Chen Pan\footnote{Department of Mathematics, USTC, P.R. China, \texttt{panchen@mail.ustc.edu.cn}}
\ and Shige Peng\footnote{School of Mathematics,
Shandong University, P.R. China
\texttt{peng@sdu.edu.cn}}}}

\maketitle

\begin{abstract}\noindent

We formulate and solve the martingale problem in a nonlinear expectation space.
 Unlike the classical work of Stroock and Varadhan (1969) where the linear operator in the associated PDE is naturally defined from the corresponding diffusion process, the main difficulty in the nonlinear setting is to identify an appropriate class of nonlinear operators for the associated fully nonlinear PDEs. 

Based on the analysis of the martingale problem,  we  introduce the notion of weak solution for stochastic differential equations under nonlinear expectations and obtain an existence theorem under  the H\"older continuity condition of the coefficients. The approach to establish the existence of weak solutions generalizes the classical Girsanov transformation method in that  it no longer requires the two (probability) measures to be absolutely continuous.


{\em MSC2000 subject classification\/}: 60G40, 60H30, 49J10, 49K10,
93E20.

\end{abstract}

\section{Introduction}\label{sec:Intro}

\subsection{Background}

A probability measure and its associated linear expectation is a
special case of nonlinear expectations.  A particular nonlinear
expectation is the \emph{sublinear} or \emph{$G$-expectation}
introduced in \cite{P07a},  defined as the following.
  Let ${{\Omega}} = C[0,\infty)$ be the space of all real valued continuous functions $({{\omega}}(t))_{t\ge0}$, and  let ${\mathcal H}$  be the linear space of random variables on ${{\Omega}}$.
A sublinear expectation $\widehat{\mathbb{E}}$ is a functional on ${\mathcal H}$ satisfying, for all $X, Y \in {\mathcal H}$,

\noindent
{\rm (S1)} Monotonicity: $\widehat{\mathbb{E}}[X] \ge \widehat{\mathbb{E}}[Y]$, if $X \ge Y$;\\
{\rm (S2)} Constant preserving: $\widehat{\mathbb{E}}[c] = c$, for $c\in {\mathbb R}$;\\
{\rm (S3)} Sub-additivity: $\widehat{\mathbb{E}}[X + Y]  \le \widehat{\mathbb{E}}[X] + \widehat{\mathbb{E}}[Y]$;\\
{\rm (S4)} Positive homogeneity: $\widehat{\mathbb{E}}[aX] = a\widehat{\mathbb{E}}[X]$, for $a \ge 0$.\\
The triple $({{\Omega}},{\mathcal H},\widehat{\mathbb{E}})$ is called a {\it sublinear expectation space}. In a sublinear
expectation space, there is no longer one-to-one correspondence
between the nonlinear expectation and its induced capacity, unlike the linear expectation and its induced probability measure.

One motivation for developing the $G$-theory is the theory
of risk measure. A coherent risk measure $\rho$ is first introduced in \cite{ADEH99}, which can be associated
with a sublinear expectation $\widehat{\mathbb{E}}$ via
$\rho(X)=\widehat{\mathbb{E}}[-X]$ for any random variable $X$.
$G$-theory provides a rigorous mathematical framework for
time-consistent risk measures, which were previously restricted to be
static. Sublinear
expectation is also related to model uncertainty. An insightful result
in \cite{DHP11} shows that a sublinear expectation is
connected to classical expectation through a class of probability
measures that measure the ``size'' of uncertainty in the following
way: there exists a weakly compact family ${\mathcal P}$ of probability measures
on $({{\Omega}}, {\cal B}({{\Omega}}))$ such that $\widehat{\mathbb{E}}[X] =
\max_{P\in{\mathcal P}} E^P[X]$ for a proper class of $X$.
 Here $\mathcal{B}({{\Omega}})$ is the Borel ${\sigma}$-algebra on ${{\Omega}}$,
and $E^P[X]$ is the linear expectation with
respect to $P$. Consequently, the notion of ``quasi-sure'' in a sublinear
expectation replaces that of ``almost-sure'' in a probability space.
From this perspective, a sublinear expectation ``measures'' the model
uncertainty: the bigger the expectation $\widehat{\mathbb{E}}$, the more
the uncertainty.

The very first building block of $G$-theory is
the {\it $G$-normal distribution}, i.e., a normal distribution with an
uncertain variance written as $N({0}\times[\underline{\sigma}^2,
\overline{\sigma}^2])$. It is characterized by the {\it $G$-heat equation}
\begin{equation}
\partial_t u -G(D^2 u)=0, \ \ u|_{t=0}=\varphi.
\end{equation}
Here, $G: {\mathbb R}\to {\mathbb R}$ is a monotonic sublinear function given by
$G(\gamma) = {\frac{1}{2}} \sup_{{\alpha} \in
  [\underline{\sigma}^2,\overline{\sigma}^2]}{\gamma {\alpha}},$ where
$\underline{\sigma}^2 = -\widehat{\mathbb{E}}[-X^2]$ and $\overline{\sigma}^2 =
\widehat{\mathbb{E}}[X^2].$  The $G$-theory is then developed in a way
similar to the classical probability theory: the notion of
$G$-(in)dependence and a $G$-central limit theorem are developed.
Especially, in order to define the conditional expectation, a backward
recursive procedure is adopted to first define a pre-expectation,
starting from the solution of the $G$-heat equation with $\varphi$.  This
idea is analogous to defining stochastic processes from a
finite-dimensional distribution.  Such a procedure is well-defined
once Kolmogorov's time-consistency theorem or the semi-group property
is established, as shown in \cite{P05}.
From here, the $G$-Brownian motion, $G$-It\^o's calculus, $G$-SDEs,
and $G$-martingale are developed similarly as the classical It\^o's
calculus. This is the $G$-theory in the spirit of Kolmogorov and It\^o.

\subsection{The martingale  problem with ${\widetilde{G}}$}
In this paper, we consider the martingale problem in the spirit of Stroock and Varadhan \cite{SV69}, albeit in a nonlinear expectation space.

\paragraph{Problem formulation}
The classical martingale problem studies
a diffusion process and its distributions with a parabolic PDE with a
linear differential operator $L_{\theta}$ and their semi-group properties,
and shows the equivalence of solving the martingale problem to the
unique weak solution of an associated stochastic differential equation
with given drift and diffusion coefficients. Moreover, the
probability measure is built along with the underlying random
processes and its uniqueness is established.  Naturally, under a nonlinear
expectation, the corresponding martingale problem is to find a family of nonlinear operators $\{{\widetilde{\mathcal E}}_t\}_{t\ge0}$ on a nonlinear expectation space $({{\Omega}}, {\mathcal H})$ such that
\begin{eqnarray} \label{equation1}
\varphi(X_t) - \int_0^t {\widetilde{G}}(X_{\theta},D\varphi(X_{\theta}),D^2\varphi(X_{\theta}))\,d\theta, t \ge 0
\end{eqnarray}
is an ${\widetilde{\mathcal E}}$-martingale for all $\varphi \in C_0^\infty({\mathbb R}^d)$. Here $X_t({\omega}) = {\omega}(t), {\omega} \in {{\Omega}}$, ${\widetilde{G}}:{\mathbb R}^d\times{\mathbb R}^d\times\mathbb{S}^d\to {\mathbb R}$ is a given continuous function satisfying certain properties to be specified later, where $\mathbb{S}^d$  is the collection of $d\times d$ symmetric matrices with usual order.

\paragraph{Appropriate class for ${\widetilde{G}}$}
However, there is a major issue.  Unlike the classical martingale
problem where the linear differential operator $L_{\theta}$ is defined
naturally as the infinitesimal generator of a diffusion process with
given drift terms $b^i(\cdot, \cdot)$ and diffusion terms $a^{ij}(\cdot, \cdot)$, the specification of
the continuous function ${\widetilde{G}}$ is not so obvious in a nonlinear setting.  Given the nonlinear
nature of the PDE, identifying the appropriate class of ${\widetilde{G}}$ is
critical for the scope and feasibility of our study.

To study the martingale problem in the nonlinear expectation space, one would need to analyze the  class of  fully nonlinear PDEs of the following form
\begin{eqnarray}
\partial_t u - {\widetilde{G}}(x,Du,D^2u) = 0,\; (t,x)\in (0,\infty)\times {\mathbb R}^d. \label{nonlinearPDE}
\end{eqnarray}
Assuming that ${\widetilde{G}}$ were chosen so that this fully nonlinear PDE had a unique solution in some sense, then one could proceed to construct a sequence of conditional expectations ${\widetilde{\mathcal E}}_t$. Now, in order for such a nonlinear expectation to be consistent with the existing literature in both linear and sublinear spaces,  it would be natural to require that ${\widetilde{G}}$  be monotonic, subadditive, and positive homogeneous. However, this class of ${\widetilde{G}}$ would appear too restrictive.
In this paper, we will show  that an appropriate class of ${\widetilde{G}}$ needs not be sublinear itself. Instead ${\widetilde{G}}$  should  be ``dominated'' by some sublinear and continuous function $G$. We call such a condition the ``DOM'' condition and
define a class $\mathfrak{D}$ for ${\widetilde{G}}$ (see Section \ref{sec:G}).

\paragraph{Fully nonlinear PDEs and weak solutions of $G$-SDEs} Once such a class of ${\widetilde{G}}$ is fixed,  then one needs to analyze the associated fully nonlinear PDE. For the pair of $G$ and ${\widetilde{G}}$, there is a pair of the associated PDEs
\eqref{PDE} and \eqref{GE2} respectively (see Section \ref{sec:pdes}).
The specification of the class $\mathfrak{D}$ as outlined in Section \ref{sec:G} has several implications. First, it ensures the uniqueness of the solutions for the PDEs. Second, it guarantees that a conditional expectation, i.e., a family of operators $\{{\widetilde{\mathcal E}}_t\}_{t\ge0}$ can be constructed
from the viscosity solution of this PDE, and the constructed conditional expectation has reasonable properties such as  time consistency. Finally,  the condition ``DOM'' not only allows one to build a piecewise Brownian motion in the ${\widetilde{G}}$ space based on the $G$-Brownian motion embedded in the sub-linear expectation space, but also ensures a much simpler stochastic calculation within the framework of Brownian motion.

Once the martingale problem is solved, it is natural to introduce the notion of weak solution of $G$-SDE under the nonlinear expectation and discuss its existence as in the classical probability theory.

\paragraph{Our work vs. related work}
The literature on sublinear expectation is growing rapidly (see \cite{P10a} and references therein).
In contrast to the bottom-up
approach in $G$-theory, where the $G$-martingale and $G$-It\^o's
calculus are developed on a sublinear expectation space from the basic $G$-heat equation, our approach starts with a general class of fully  nonlinear PDEs which includes the $G$-heat equation as a very special case.  These PDEs are state dependent. Consequently, our analysis on the existence and uniqueness of their solutions not only generalizes existing  results  including  \cite{P10a} and \cite{FS92}, but also leads to the construction of nonlinear expectations that goes beyond the sublinear ones in \cite {P10a}, \cite{STZ12} and \cite{DHP11}. As a result, random processes, especially  martingales, and the stochastic calculus are all developed under a more general framework.

It is worth mentioning that there are other approaches in addition to ours of constructing nonlinear expectations from PDEs.  For example, \cite{STZ12} uses the classical stochastic control approach of building regular conditional probability when considering a family of backward stochastic differential equations; and \cite{N12} takes a  control approach, where random $G$-expectations are constructed based on an optimal control formulation with path-dependent control sets, hence a path-dependent family of probabilities. All these approaches, however, leads to sublinear expectation spaces (see also \cite{DHP11}).

Finally, our approach to establish the existence of weak solutions generalizes the classical Girsanov transformation method in that  it no longer requires the two (probability) measures to be absolutely continuous. Instead it is critical  that probability measures singular from each other should be ``dominated'' by a certain sublinear expectation.

\paragraph{Outline of the paper}
The paper starts with the discussion of $G$ and ${\widetilde{G}}$ in the class  $\mathfrak{D}$ (Section \ref{sec:G}) and states some of the key properties of the solutions to the associated PDEs (Section \ref{sec:pdes}). Following the Kolmogorov's idea,  a family of nonlinear operators $\{{\widetilde{\mathcal E}}_t\}_{t\ge0}$ is constructed via solutions of the PDEs \eqref{GE2} and \eqref{PDE} with analysis of their properties (Section \ref{sec:NE}); Section~\ref{sec:mp} finishes the proof of the martingale problems. Finally, the weak solution of $G$-SDE is introduced and analyzed in Section \ref{sec:weak}. Appendix contains some technical details for stochastic calculus under nonlinear expectations and proofs on the existence and uniqueness for the PDEs that are associated with the martingale problem.

\section{Martingale problems}
\label{sec:MP}

\subsection{$G$ and ${\widetilde{G}}$}
\label{sec:G}

Let us start by defining a class $\mathfrak{D}$ of functions ${\widetilde{G}}$. ${\widetilde{G}}$ is continuous and``dominated'' by a continuous and sublinear function $G$.

\begin{definition}[Class $\mathfrak{D}$]
A continuous function ${\widetilde{G}}:{\mathbb R}^d\times{\mathbb R}^d\times\mathbb{S}^d\to {\mathbb R}$ is of {\em class $\mathfrak{D}$}, if there exists a  continuous function $G: {\mathbb R}^d\times{\mathbb R}^d\times\mathbb{S}^d\to {\mathbb R}$ such that
\begin{itemize}
 \item ${\widetilde{G}}(x,{\alpha} p,{\alpha} A)= {\alpha} {\widetilde{G}}(x,p,A)$ for all $x\in{\mathbb R}^d, {\alpha} \ge 0$,
 \item  for each $x\in{\mathbb R}^d$, $(p,A),(p',A')\in {\mathbb R}^d\times\mathbb{S}^d$,
 \begin{equation*}
\widetilde{G}(x,p,A) - \widetilde{G}(x,p',A') \le G(x,p-p',A-A'), \tag{DOM}\label{DOM}
\end{equation*}
\end{itemize}
with $G$ satisfying
\begin{enumerate}[(A).]
  \item  (Subadditivity) $G(x,p+\bar{p},A+\bar{A}) \le G(x,p,A) + G(x,\bar{p},\bar{A})$, \label{c1}
  \item  (Positive Homogeneity) $G(x, \beta p, \beta A) = \beta G(x,p,A), \beta \ge 0$, \label{c2}
  \item (Monotonicity) $G(x,p,A) \le G(x,p,A+\bar{\bar{A}})$, \label{c3}
  \item (Uniform  Lipschitz Continuity) $|G(x,p,A) - G(x,\bar{p},\bar{A})| \le L(|p-\bar{p}|+|A-\bar{A}|)$,  for some  $L >0$, \label{c4}
\end{enumerate}
for any $ x,p,\bar{p}\in {\mathbb R}^d$, $A,\bar{A},\bar{\bar{A}} \in \mathbb{S}^d$, and $\bar{\bar{A}}\ge0$.
\end{definition}

\begin{example} \label{example1}
 When $G$ satisfies conditions \eqref{c1} and \eqref{c2}, then by Theorem I.2.1 of ~\cite{P10a}, for each given $x\in{\mathbb R}^d$, there exists a bounded, closed, and convex subset $U(x)\subset\mathbb{S}^d\times{\mathbb R}^d$, such that
\begin{equation*}
G(x,p,A) = \sup_{(a,b)\in U(x)}\left\{{\frac{1}{2}}\text{\em tr}[aA] + \langle b,p\rangle\right\}.
\end{equation*}
Since for each given $x$, there exists a dense subset $U_0(x)\subset U(x)$, which is countable, also denoted by $\left\{({\frac{1}{2}} a(x,i),b(x,j))\right\}_{i,j\in{\mathbb{N}}}$, one can rewrite the above expression in the following sublinear form
\begin{eqnarray}\label{expofG}
G(x,p,A) = \sup_{\gamma\in \Gamma}\left\{{\frac{1}{2}}\text{\em tr}[a(x,\gamma)A] + \langle b(x,\gamma),p\rangle\right\},
\end{eqnarray}
where $\Gamma$ is an index set, and $a(x,\gamma)\in\mathbb{S}^d,b(x,\gamma)\in{\mathbb R}^d$ are bounded. Moreover, when $G$ also satisfies \eqref{c3}, then $a(x,\gamma)\ge0$ for any $x\in{\mathbb R}^d$ and $\gamma\in\Gamma$. And there exists $\sigma(x,\gamma)\in{\mathbb R}^{d\times d}$, such that $a(x,\gamma) = {\sigma}(x,\gamma){\sigma}^T(x,\gamma).$
With additional condition \eqref{c4}, one can simply write
 \begin{eqnarray}\label{expofG1}
G(p,A) = \sup_{(\gamma,\beta)\in{\Gamma}}\left\{{\frac{1}{2}} \text{\em tr}[{\gamma}A]+\langle \beta, p \rangle\right\}
\end{eqnarray}
with ${\Gamma}\subset \mathbb{S}^d_+\times{\mathbb R}^d$ being bounded, convex, and closed.

Of course, such a $G$ is dominated by itself in the sense of \eqref{DOM}, thus  in the class $\mathfrak{D}$.
\end{example}

\begin{example}\label{exoftG}
Assume ${\widetilde{G}}$ of the form
  \begin{align}
  {\widetilde{G}}(x,p,A) = \sup_{{\gamma}\in {\Gamma}}\inf_{{\lambda}\in {\Lambda}}\left\{{\frac{1}{2}}\text{\em tr}[{\sigma}(x,\gamma,{\lambda}){\sigma}^T(x,\gamma,{\lambda})A]+\langle b(x,\gamma,{\lambda}),p\rangle \right\} \label{form1}
  \intertext{or}
  {\widetilde{G}}(x,p,A) = \inf_{{\gamma}\in{\Gamma}}\sup_{{\lambda}\in {{\Lambda}}}\left\{{\frac{1}{2}}\text{\em tr}[{\sigma}(x,\gamma,{\lambda}){\sigma}^T(x,\gamma,{\lambda})A]+\langle b(x,\gamma,{\lambda}),p\rangle \right\}.\label{form2}
  \end{align}
   Here ${\Gamma}$ and ${{\Lambda}}$ are compact metric spaces, ${\sigma},b\in C_b({\mathbb R}^d\times{\Gamma}\times{{\Lambda}})$, and ${\sigma}(\cdot, {\gamma},{\lambda})$ and $b(\cdot,{\gamma},{\lambda})$ are uniformly Lipschitz continuous in the following sense
  \begin{equation*}
  |{\sigma}(x,\gamma,{\lambda}) - {\sigma}(y,{\gamma},{\lambda})| + |b(x,{\gamma},{\lambda}) - b(y,{\gamma},{\lambda})| \le \tilde{L}|x-y|, \text{ for all ${\gamma}\in{\Gamma}, {\lambda} \in{\Lambda}$},
  \end{equation*}
  with $\tilde{L}>0$ a constant.

  Such form of ${\widetilde{G}}$ is important for  stochastic games \cite{BBP97} and is in class  $\mathfrak{D}$:  it is clearly dominated by $G$ specified by
  \begin{equation*}
  G(x,p,A) = \sup_{{\gamma}\in{\Gamma},{\lambda}\in{\Lambda}} \left\{ {\frac{1}{2}}\text{\em tr}[{\sigma}(x,{\gamma},{\lambda}){\sigma}^T(x,{\gamma},{\lambda})A]+\langle b(x,\gamma,{\lambda}),p\rangle \right\}.
  \end{equation*}
\end{example}

Throughout the paper unless otherwise specified, $G$ and ${\widetilde{G}}$ satisfy the conditions specified in the definition of class  $\mathfrak{D}$. And without loss of generality, we assume $G$ is of the form (\ref{expofG1})
as discussed in Example~\ref{example1}.

\subsection{PDEs associated with $G$ and ${\widetilde{G}}$}
\label{sec:pdes}

Now we introduce two classes of fully nonlinear PDEs associated with $G$ and ${\widetilde{G}}$ in class $\mathfrak{D}$.

\paragraph {State-dependent parabolic PDEs associated with $G$}

\begin{equation*}
\begin{cases}\label{PDE}
\partial_t u(t,x) - G(x,Du(t,x),D^2u(t,x)) = 0,\; (t,x)\in (0,\infty)\times {\mathbb R}^d,\\
u(0,x) = \varphi(x), \;x\in{\mathbb R}^d.
\end{cases}\tag{P}
\end{equation*}
\paragraph {Fully nonlinear PDEs associated with ${\widetilde{G}}$}
\begin{equation*}
\begin{cases}
  \partial_t u(t,x) - \widetilde{G}(x,Du(t,x),D^2u(t,x)) = 0, (t,x)\in(0,\infty)\times{\mathbb R}^d,\\
  u(0,x) = \varphi(x), x\in {\mathbb R}^d.
\end{cases}\tag{$\widetilde{\text{P}}$}\label{GE2}
\end{equation*}

Clearly the PDE ~\eqref{PDE} is a special case of the PDE ~\eqref{GE2}.

The conditions specified in class $\mathfrak{D}$ for $G$ and ${\widetilde{G}}$ are essential to ensure the existence of the viscosity solutions for PDE ~\eqref{PDE} and PDE ~\eqref{GE2}  and to guarantee that such solutions have  nice properties.

In the following, we will discuss the existence, the uniqueness, and the properties of the solutions associated with   PDEs ~\eqref{PDE} and~\eqref{GE2}. Note that some results hold under more general conditions for $G$ and ${\widetilde{G}}$ than those specified in class $\mathfrak{D}$. To avoid confusions, all conditions in the theorems and lemmas are specified for $G$ and ${\widetilde{G}}$.

First, recall the definition of the viscosity solutions for  the associated PDEs in~\eqref{GE2} and~\eqref{PDE}.

\begin{definition}
  Given a constant $T>0$. A {\em viscosity subsolution} of the PDE in \eqref{GE2} on $(0,T)\times{\mathbb R}^d$ is an upper semicontinuous  (USC) function $u$ in $(0,T)\times{\mathbb R}^d$ such that for all $(t,x)\in (0,T)\times {\mathbb R}^d,\phi \in C^2((0,T)\times{\mathbb R}^d)$ such that $u(t,x) = \phi(t,x)$ and $u<\phi$ on $(0,T)\times{\mathbb R}^d\setminus\{(t,x)\}$, we have
 \begin{equation*}
  \partial_t \phi(t,x) - {\widetilde{G}}(x,D\phi(t,x),D^2\phi(t,x)) \le 0;
  \end{equation*}
  likewise, a {\em viscosity supersolution} of  the PDE in \eqref{GE2} on $(0,T)\times{\mathbb R}^d$ is a lower semiconitouns (LSC) function  $v$ in $(0,T)\times{\mathbb R}^d$ such that for all $(t,x)\in (0,T)\times {\mathbb R}^d,\psi \in C^2((0,T)\times{\mathbb R}^d)$ such that $v(t,x) = \psi(t,x)$ and $v>\psi$ on $(0,T)\times{\mathbb R}^d\setminus\{(t,x)\}$, we have
  \begin{equation*}
  \partial_t \psi(t,x) - {\widetilde{G}}(x,D\psi(t,x),D^2\psi(t,x)) \ge 0.
  \end{equation*}
  And a {\em viscosity solution} of  the PDE in \eqref{GE2} on $(0,T)\times{\mathbb R}^d$ is a function that is both a viscosity subsolution and a viscosity supersolution of  the PDE in \eqref{GE2} on $(0,T)\times{\mathbb R}^d$.
\end{definition}
The definition of the viscosity solution to PDE in~\eqref{PDE} is similar, with ${\widetilde{G}}$ replaced by $G$.

Now, note that PDE \eqref{PDE} has been extensively studied, for example, in the literature of portfolio selections
  (see for instance \cite{Ph09}). And its comparison theorem can be established  with slightly modified techniques from \cite{FS92}.
\begin{theorem}[Comparison theorem for PDE in \eqref{PDE}]\label{Comparison}
  Given a continuous function $G: {\mathbb R}^d\times{\mathbb R}^d\times\mathbb{S}^d\to {\mathbb R}$, which satisfies conditions \eqref{c1}, \eqref{c2}, \eqref{c3}. Suppose ${\sigma},b$ are uniformly Lipschitz continuous with respect to $x$.
  Let $\underline{u} \in \text{USC}([0,T]\times{\mathbb R}^d)$ be a viscosity subsolution of the PDE in \eqref{PDE} and $\overline{u} \in \text{LSC}([0,T]\times{\mathbb R}^d)$ be a viscosity supersolution of the PDE in \eqref{PDE} on $[0,T]\times{\mathbb R}^d$ with polynomial growth. Then $\underline{u} \le \overline{u}$ when $\underline{u}|_{t=0} \le \overline{u}|_{t=0}$.
\end{theorem}


Perron's existence result of the solution of \eqref{PDE} follows from  Appendix C.3 of \cite{P10a}.  

\begin{theorem}[Existence for PDE \eqref{PDE}]\label{existofsol}
  Assuming a comparison theorem holds for \eqref{PDE}. Moreover, suppose that there is a viscosity subsolution of \eqref{PDE} $\underline{u}$ and a viscosity supersolution $\overline{u}$ of \eqref{PDE} such that  $\underline{u}_*|_{t=0} = \overline{u}^*|_{t=0} = \varphi \in C({\mathbb R}^d)$ with polynomial growth.  Here $u^*$ is  the upper semicontinuous envelope of $u$ and $u_*$ is lower semicontinuous envelope of $u$.
  Then
  \[
  w(t,x) = \sup\{W(t,x);\underline{u}\le W \le \overline{u} \text{ and $W$ is a viscosity subsolution of \eqref{PDE}}\},
  \]
  is a viscosity solution of \eqref{PDE}.

  In particular, if problem~\eqref{PDE} satisfies conditions \eqref{c1}, \eqref{c2}, and \eqref{c3}, and ${\sigma},b$ are bounded and uniformly Lipschitz continuous, then it has a unique solution.
\end{theorem}

We next state the comparison theorem for PDE ~\eqref{GE2} which relies on a technical condition ~\eqref{homo} as detailed in the Appendix. We also outline its proof in the Appendix.

\begin{theorem}[Comparison theorem for PDE in~\eqref{GE2}] \label{comparison3}
  Suppose both $G$ and ${\widetilde{G}}$ satisfy condition ~\eqref{homo}, with  their respective corresponding continuous decomposition functions satisfying condition~\eqref{g}. Suppose  ${\widetilde{G}}$ satisfies condition~\eqref{DOM} and $G$ satisfies conditions \eqref{c3} and \eqref{c4}. Let $\underline{\widetilde{u}}\in USC([0,T]\times{\mathbb R}^d)$ be a subsolution of the PDE in \eqref{GE2} and $\overline{\widetilde{u}}\in LSC([0,T]\times{\mathbb R}^d)$ be a supersolution of the PDE in \eqref{GE2} on $(0,T)\times{\mathbb R}^d$ and $\overline{w}$ is a supersolution of the PDE in \eqref{PDE}. They all satisfy the polynomial growth condition. If $(\underline{\widetilde{u}}-\overline{\widetilde{u}})|_{t=0} = \overline{w}|_{t=0}$, then $\underline{\widetilde{u}}-\overline{\widetilde{u}} \le  \overline{w}$ on $[0,T)\times{\mathbb R}^d$.

  In particular, $\underline{\widetilde{u}}\le\overline{\widetilde{u}}$ on $[0,T)\times{\mathbb R}^d$ provided that $\underline{\widetilde{u}}|_{t=0}\le \overline{\widetilde{u}}|_{t=0}$.
\end{theorem}

Moreover,  the same proof of Theorem C.3.1 of~\cite{P10a} leads to
\begin{theorem}[Existence of the solution of \eqref{GE2}]
Suppose $G$ satisfies conditions~\eqref{c1}, \eqref{c2}, \eqref{c3}, and \eqref{c4}. Assume that both $G$ and ${\widetilde{G}}$ satisfy condition~\eqref{homo}, with their respective corresponding continuous decomposition functions satisfying condition~\eqref{g}. If  ${\widetilde{G}}$ satisfies the ~\eqref{DOM} condition, the relation
\[
\widetilde{G}(x,\alpha p,\alpha A) - \widetilde{G}(x,0,0) = \alpha[\widetilde{G}(x, p, A)-\widetilde{G}(x,0,0)]\quad \text{for all $\alpha \ge 0$},
\]
and $\widetilde{G}(x,0,0)$ has polynomial growth, then there exists a unique solution for PDE ~\eqref{GE2}.
\end{theorem}


\begin{remark}
   In fact, the positive homogeneity condition and condition \eqref{homo} are not necessary for the uniqueness of the solution for PDE ~\eqref{GE2}. For instance, take
  \begin{eqnarray}\label{formoftG2}
  {\widetilde{G}}(x,p,A) = \sup_{\gamma\in\Gamma} \{a({\gamma})A + g(x,p,{\gamma})\}, (x,p,A) \in {\mathbb R}\times{\mathbb R}\times{\mathbb R},
  \end{eqnarray}
  where $\Gamma$ is an index set such that $a(\gamma)\ge 0$ is uniformly bounded, and the continuous function $g$ is dominated by a continuous function $h:{\mathbb R}\times{\mathbb R}\to {\mathbb R}$ in the  sense of
 \begin{equation*}
  g(x,p,{\gamma}) - g(x,\bar{p},{\gamma}) \le h(x,p-\bar{p}), \quad\text{for every ${\gamma}\in{\Gamma}$}.
  \end{equation*}
  Suppose $g$ satisfies $|g(x,p,{\gamma})-g(y,p,{\gamma})| \le L_g(1+|x|+|y|)|x-y|(1+|p|), x,y,p\in{\mathbb R}$ uniformly in ${\gamma}\in{\Gamma}$ and $L_g>0$ is a constant, and $|h(x,p)-h(y,p)| \le L_h(1+|x|+|y|)|x-y|(1+|p|)$.
  Then this ${\widetilde{G}}$ is dominated by $G(x,p,A) = \sup_{\gamma\in\Gamma} \{a({\gamma})A\} + h(x,p)$, yet ${\widetilde{G}}$ does not satisfy the positive homogeneity condition and condition~\eqref{homo}.
  Nevertheless, the same approach shows that the comparison still holds.
\end{remark}


For ${\widetilde{G}}$ of the forms as in Example~\ref{exoftG}, we have the following results about the associated PDEs.

\begin{theorem}\label{domthm}\cite{BBP97}
  Suppose ${\widetilde{G}}$ is of the form \eqref{form1}.  Let $u_1, -u_2 \in \text{USC}([0,T]\times{\mathbb R}^d)$, $u_1$ and $u_2$ are subsolution and supersolution to problem~\eqref{GE2} with initial conditions $u_1(0,x) = \varphi_1(x)$ and $u_2(0,x) = \varphi_2(x)$, respectively. Set $w = u_1 - u_2$. Then $w$ is a subsolution to
  \begin{eqnarray}
  \begin{cases}
  \partial_t w(t,x) - G(x,Dw(t,x),D^2w(t,x)) = 0, \, (t,x)\in (0,T]\times{\mathbb R}^d,\\
  w(0,x) = \varphi_1(x) - \varphi_2(x),\, x\in{\mathbb R}^d,
  \end{cases}
  \end{eqnarray}
  with $\displaystyle G(x,p,A) = \sup_{{\gamma} \in {\Gamma}, {\lambda} \in {\Lambda}} \left\{ {\frac{1}{2}} \text{\em tr}[{\sigma}(x,{\gamma},{\lambda}){\sigma}^T(x,{\gamma},{\lambda})A] + \langle b(x,{\gamma},{\lambda}),p\rangle \right\}$.
  Consequently, let $u\in \text{\em USC}([0,T]\times{\mathbb R}^d), v\in \text{\em LSC}([0,T]\times{\mathbb R}^d)$ be the subsolution and supersolution to the PDE~\eqref{GE2}, respectively. If both $u$ and $v$ have at most polynomial  growth and $u(0,x)\le v(0,x)$, then $u(t,x)\le v(t,x)$ in $(0,T]\times{\mathbb R}^d$.
\end{theorem}
 Since our PDE has a simpler form than that in \cite{BBP97} without the jump term, their proof can be greatly simplified, as illustrated in the Appendix.  Similar results and proof hold for PDEs with ${\widetilde{G}}$ of form \eqref{form2}.

Finally, we discuss  the properties of the solutions for the PDEs, when exist. Clearly, from the definition of class $\mathfrak{D}$, one sees

\begin{theorem}[Properties of the solutions of PDEs]\label{PropofSol}
 Let $u^\phi,\tilde{u}^\varphi \in C([0,T]\times{\mathbb R}^d)$ denote the unique solutions of \eqref{PDE} and \eqref{GE2} with polynomial growth, with the boundary conditions $\phi$ and $\varphi$ respectively. Then we have
  \begin{align*}
    \tilde{u}^{\varphi + c} &= \tilde{u}^\varphi + c,\\
    \tilde{u}^\varphi - \tilde{u}^\phi &\le u^{\varphi-\phi},\\
    \tilde{u}^{{\alpha} \varphi} &= {\alpha} \tilde{u}^\varphi,\: {\alpha} \ge 0,
  \end{align*}
  where $c\in{\mathbb R}$ is a constant, and $\varphi, \phi$ are continuous functions with polynomial growth.
\end{theorem}




\subsection{Nonlinear expectations ${\widetilde{\mathcal E}}$ and ${\mathcal E}$}
\label{sec:NE}
\subsubsection{Construction of ${\widetilde{\mathcal E}}$ and ${\mathcal E}$ from the associated PDEs}
\label{construction}

Assuming the unique solution $\tilde{u}^\varphi$ to PDE ~\eqref{GE2},  one can define the `conditional expectation' ${\widetilde{\mathcal E}}_t$ for $t\in [0, T]$, $T<\infty$.

The construction starts from the ``pre-expectation''. Let ${\Omega} = C_x([0,\infty);{\mathbb R}^d) = \{{\omega}(\cdot);{\omega} \text{ is a}$ continuous ${\mathbb R}^d$-valued function on $[0,\infty)$ and ${\omega}(0)=x\}$. Fix $N>0$, take $0 = t_0\le t_1 \le \cdots \le t_N \le T.$
 Take  $\varphi_0$ from  a proper function space on $({\mathbb R}^d)^N$ denoted by $\mathcal{C}(({\mathbb R}^d)^N)$, and set $\xi({\omega})=\varphi_{0}(X_{t_{1}},\cdots,X_{t_{N}})$. Let ${\mathcal T}_t[\varphi(\cdot)](x) := u(t,x)$, and
   for $0\le j\le N$, define
\begin{align*}
\varphi_{1}(x_{1},\cdots,x_{N-1}) &  ={\mathcal T}_{t_N - t_{N-1}}[\varphi_{0}(x_{1},\cdots
,x_{N-1},\cdot)](x_{N-1}),\\
&  \vdots\\
\varphi_{N-j}(x_{1},\cdots,x_{j}) &  ={\mathcal T}_{t_{j+1}-t_j}[\varphi_{N-j-1}(x_{1},\cdots
,x_{j},\cdot)](x_{j}),\\
&  \vdots\\
\varphi_{N-1}(x_{1}) &  ={\mathcal T}_{t_2 - t_1}[\varphi_{N-2}(x_{1},\cdot)](x_{1}),\\
\varphi_{N} &  ={\mathcal T}_{t_1}[\varphi_{N-1}(\cdot)](x),
\end{align*}
where $\varphi_k \in \mathcal{C}(({\mathbb R}^d)^{N-k}),\; 0 \le k \le N-1$ and $\varphi_N \in {\mathbb R}$.
Then define
\begin{equation*}
{\widetilde{\mathcal E}}_t[\xi] = \varphi_{N-j}(X_{t_1},\cdots,X_{t_j}), \text{ if $t=t_j,0\le j\le N.$}
\end{equation*}

This construction approach is in the spirit of Nisio's semigroup theory (see \cite{N76a}, \cite{N76b}, and \cite{P05}), where nonlinear expectations are constructed from  nonlinear Markov chains after establishing a generalized Kolmogorov's consistency theorem and pre-expectations.

In our case, we denote such an `expectation' by ${\widetilde{\mathcal E}}$.
Since PDE ~\eqref{PDE} is a special case of PDE ~\eqref{GE2},  a sublinear expectation ${\mathcal E}$ can be similarly defined from the its solution. It is a special case of ${\widetilde{\mathcal E}}$.

\begin{remark}
\begin{enumerate}[1).]
  \item One can set $\mathcal{C}(({\mathbb R}^d)^N) := C_{l.Lip}(({\mathbb R}^d)^N)$, where $C_{l.Lip}({\mathbb R}^n)$ is the space of real valued continuous functions defined on ${\mathbb R}^n$ such that
\begin{equation*}
|\varphi(x) - \varphi(y)| \le C(1 + |x|^m + |y|^m)|x - y|, \forall \, x, y \in {\mathbb R}^n,
\end{equation*}
for some $C > 0$, and $m \in {\mathbb{N}}$ depending on $\varphi$.

 Now let ${\Omega}_T:= \{{\omega}({\cdot\wedge T});{\omega}\in{\Omega}\}$ and ${\mathcal H} := L_{ip}({\Omega}_T) = \{\varphi(X_{t_{1}},\cdots,X_{t_{N}});\varphi\in C_{l.Lip}(({\mathbb R}^d)^N)$ for some $N\in{\mathbb{N}} \text{ and}\; 0 \le t_1 \le \cdots \le t_N \le T\}$, and call $({\Omega},{\mathcal H},{\widetilde{\mathcal E}})$ a nonlinear expectation space. If $\xi = \varphi(X_{t_{1}},\cdots,X_{t_{N}})$ with $t_N\le t \in[0,T]$, we say $\xi\in L_{ip}({\Omega}_t)$. It is clear that $L_{ip}({\Omega}_t)\subset L_{ip}({\Omega}_T), t\le T$.

\item With a sublinear ${\mathcal E}$, for each $t\in[0,T]$, one can extend the space $L_{ip}({\Omega}_t)$ to a Banach space $L^1_{\mathcal E}({\Omega}_t)$ under the norm $\|\cdot\| : ={\mathcal E}[|\cdot|]$ as  in~\cite{P10a}, or see the Appendix, since the nonlinear expectation ${\mathcal E}$ is sublinear. And from now on, we take ${\mathcal H} = L^1_{\mathcal{E}}({\Omega}):= \cup_{T>0}L^1_{\mathcal{E}}({\Omega}_T)$.

\end{enumerate}
\end{remark}

\begin{remark}\label{cexam}
   We point out here the above construction procedure for ${\widetilde{\mathcal E}}$ depends critically on the homogeneity of the PDE ~\eqref{GE2}. Otherwise, the  resulting ${\widetilde{\mathcal E}}$  may  not be well defined as seen from the following  example.
   Consider the following linear nonhomogeneous PDE
    \begin{equation*}
      \begin{cases}
        \partial_t u - \partial_{x}^2u = x,\; (t,x)\in (0,T]\times\mathbb{R},\\
        u(0,x) = c, x\in{\mathbb R}.
      \end{cases}
    \end{equation*}
    The solution is $u(t,x) = tx + c$.
   If we were to define  ${\widetilde{\mathcal E}}$ as suggested in the above procedure, and consider the constant as a function
  \begin{equation*}
  \phi_0 : \mathbb{R}^2 \to \{c\}; (x_1,x_2)\mapsto \phi_0(x_1,x_2).
  \end{equation*}
  Then
 \begin{equation*}
  c = \phi_0(X_{t_1},X_{t_2}),\; \forall t_1\le t_2 \in (0,T].
  \end{equation*}
Now clearly
  \begin{align*}
    {\widetilde{\mathcal E}}_{t_1}[c] &= \mathcal{T}_{t_2-t_1}[\phi_0(x_1,\cdot)]_{x_1=X_{t_1}}=[(t_2-t_1)x_1 + c]|_{x_1 = X_{t_1}} =: \phi_1(x_1)|_{x_1 = X_{t_1}},\\
    {\widetilde{\mathcal E}}_{0}[\phi_1(X_{t_1})] &= t_1 X_0 + (t_2 - t_1) X_0 + c = t_2 X_0 +c.
  \end{align*}
  We would have  ${\widetilde{\mathcal E}}[c] = t_2 X_0 + c$ for any arbitrary $t_2$. Thus the nonlinear expectation is not well defined.
\end{remark}
\subsubsection{Properties of nonlinear expectations ${\widetilde{\mathcal E}}$ and ${\mathcal E}$}

 Given a nonlinear expectation space $({\Omega},{\mathcal H},{\widetilde{\mathcal E}})$, a  stochastic process $(\xi_t)_{t\ge0}$ is a collection of random variables on $({\Omega},{\mathcal H})$. That is, for each $t\ge0$, $\xi_t\in {\mathcal H}$. Moreover,

\begin{definition}[${\widetilde{\mathcal E}}$-Martingale]
A stochastic process $(M_t)_{t\ge0}$ is called an {\it ${\widetilde{\mathcal E}}$-martingale} if for each $t\in[0,\infty), M_t\in L_{\mathcal E}^1({\Omega}_t)$, and for each $s\in[0,t]$,
 \begin{equation*}
  {\widetilde{\mathcal E}}_s[M_t] = M_s.
  \end{equation*}
\end{definition}

\begin{remark}
In this paper, since ${\widetilde{\mathcal E}}$  is constructed from the PDEs associated with ${\widetilde{G}}$, sometimes the ${\widetilde{\mathcal E}}$-martingale is also referred to   ${\widetilde{G}}$-martingale when there is no risk of confusion.
\end{remark}

Moreover,  the ${\widetilde{\mathcal E}}$ and ${\mathcal E}$ constructed in Section \ref{construction} have the following properties.

\begin{proposition}\label{PropofExp}
  Given a nonlinear expectation space $({\Omega},{\mathcal H},{\widetilde{\mathcal E}})$, let $\xi,\eta \in {\mathcal H}.$
  \begin{enumerate}[(I)]
    \item\label{P1} For $\varphi \in C_{l.Lip}({\mathbb R})$ and $s\le t$
 \begin{equation*}
  {\widetilde{\mathcal E}}_s[\varphi(X_t)] = u^\varphi(t-s,X_s).
  \end{equation*}

    \item\label{P2} (Monotonicity) If $\xi \le \eta$,
      \[
      {\widetilde{\mathcal E}}_t[\xi] \le {\widetilde{\mathcal E}}_t[\eta].
      \]

  \item\label{P3} (Constant preserving)
If $\xi \in L_{ip}({\Omega}_s), \eta \in
      L_{ip}({\Omega}_{s+h}), s,h \ge 0$,
      \begin{equation*}
      {\widetilde{\mathcal E}}_s[\xi + \eta] = \xi + {\widetilde{\mathcal E}}_s[\eta].
      \end{equation*}
     In particular, ${\widetilde{\mathcal E}}[\xi + c] = {\widetilde{\mathcal E}}[\xi] + c$, with $c\in{\mathbb R}$ a constant.

  \item\label{P4} (Tower property) For any $s, h >0$,
  \begin{equation*}
   {\widetilde{\mathcal E}}_s\circ{\widetilde{\mathcal E}}_{s+h} = {\widetilde{\mathcal E}}_s.
  \end{equation*}
  \item\label{domination} (Domination)
  \begin{equation*}
  {\widetilde{\mathcal E}}_t[\xi] - {\widetilde{\mathcal E}}_t[\eta] \le {\mathcal E}_t[\xi - \eta].
  \end{equation*}

  \end{enumerate}
  In addition, we have for the sublinear expectation ${\mathcal E}$,
  \begin{enumerate}[(I)]
  \setcounter{enumi}{5}
  \item\label{P6} (Subadditivity)
  \begin{equation*}
  {\mathcal E}_t[\xi + \eta] \le {\mathcal E}_t[\xi] + {\mathcal E}_t[\eta].\end{equation*}
  \item\label{P7} (Positive homogeneity)

If $\xi \in L_{ip}({\Omega}_s), \eta \in
      L_{ip}({\Omega}_{s+h}), s,h\ge 0$,
      \begin{equation*}
      {\mathcal E}_s[\xi\eta] = \xi^+{\mathcal E}_s[\eta]+\xi^-{\mathcal E}_s[-\eta].
      \end{equation*}
      In particular, for any constant $\lambda \ge 0$,
  \begin{equation*}
  {\mathcal E}_t[\lambda\xi] = \lambda{\mathcal E}_t[\xi].
  \end{equation*}
\end{enumerate}
\end{proposition}

\noindent{\bf Proof.}
For \eqref{P1}, note  that $\varphi(X_t) = \varphi(X_t+X_s-X_s) =: \psi(X_s,X_t)$,
  \begin{equation*}
    {\widetilde{\mathcal E}}_s[\varphi(X_t)] = {\widetilde{\mathcal E}}_s[\psi(X_s,X_t)] = {\mathcal T}_{t-s}[\psi(x_1,\cdot)](x_1)|_{x_1 = X_s}= u^\varphi(t-s,x_1)|_{x_1 = X_s}.
  \end{equation*}

\eqref{P2} is an implication of the comparison theorem for the PDE in \eqref{GE2}.

For \eqref{P3}, without loss of generality, assume $\xi = \varphi(X_s), \eta=\phi(X_{s+h})$, $\varphi,\phi\in\mathcal{C}({\mathbb R}^d).$ Then
\begin{equation*}
\begin{aligned}
  {\widetilde{\mathcal E}}_s[\xi + \eta] =\,& {\widetilde{\mathcal E}}_s[\varphi(X_s) + \phi(X_{s+h})]\\
   =\,& \mathcal{T}_h[\varphi(y) +\phi(\cdot)](y)|_{y=X_s}\\
   =\,& \varphi(X_s) + \mathcal{T}_h[\phi(\cdot)](y)|_{y=X_s}\\
   =\,& \xi+{\widetilde{\mathcal E}}_s[\eta],
\end{aligned}
\end{equation*}
where the third equality follows from Theorem~\ref{PropofSol}.

To show \eqref{P4}, without loss of generality, we assume $\xi = \varphi_0(X_{t_1},X_{t_2},X_{t_3}),t_1=s,t_2=s+h,t_2\le t_3\le T,\varphi_0 \in \mathcal{C}(({\mathbb R}^d)^3)$. Then
\begin{equation*}
{\widetilde{\mathcal E}}_{t_2}[\xi] = \varphi_1(X_{t_1},X_{t_2}) \text{ and} \;{\widetilde{\mathcal E}}_{t_1}[\xi] = \varphi_2(X_{t_1}), \varphi_1 \in \mathcal{C}(({\mathbb R}^d)^2),\varphi_2 \in\mathcal{C}({\mathbb R}^d),
\end{equation*}
from the construction procedure for ${\widetilde{\mathcal E}}$. Meanwhile,
\begin{equation*}
{\widetilde{\mathcal E}}_{t_1}[{\widetilde{\mathcal E}}_{t_2}[\xi]] = {\widetilde{\mathcal E}}_{t_1}[\varphi_1(X_{t_1},X_{t_2})] = \bar{\varphi}_2(X_{t_1}),
\end{equation*}
with $\bar{\varphi}_2(x) = \mathcal{T}_{t_2-t_1}[\varphi_1(x,\cdot)](x) = \varphi_2(x)$, since the PDE ~\eqref{GE2} has a unique solution for a given appropriate initial condition.

\eqref{domination} can be derived directly from Theorem~\ref{PropofSol} and the construction of ${\widetilde{\mathcal E}}$, while \eqref{P6} is a special case of \eqref{domination}.

To prove \eqref{P7}, assume $\xi = \varphi(X_s), \eta=\phi(X_{s+h})$, $\varphi,\phi\in\mathcal{C}({\mathbb R}^d).$ Then, by Theorem~\ref{PropofSol},
\begin{equation*}
\begin{aligned}
  {\mathcal E}_s[\xi\eta]&= {\mathcal E}_s[\varphi(X_s)\phi(X_{s+h})]\\
   &= \mathcal{T}_h[\varphi(y)\phi(\cdot)](y)|_{y=X_s}\\
   &= \mathcal{T}_h[\varphi(y)^+\phi(\cdot)-\varphi(y)^-\phi(\cdot)](y)|_{y=X_s}\\
   &=\{\varphi(y)^+\mathcal{T}_h[\phi(\cdot)]+\varphi(y)^-\mathcal{T}_h[-\phi(\cdot)](y)\}|_{y=X_s}\\
   &=\xi^+\mathcal{T}_h[\phi(\cdot)]_{y=X_s} + \xi^-\mathcal{T}_h[-\phi(\cdot)]_{y=X_s}\\
   &=\xi^+{\mathcal E}_s[\eta]+\xi^-{\mathcal E}_s[-\eta].
\end{aligned}
\end{equation*}
\endofproof

In addition, both ${\widetilde{\mathcal E}}$ and ${\mathcal E}$ enjoy the following property.

\begin{lemma}\label{PropI.3.7}
  Let $\xi,\eta$ be two random variables in the sublinear expectation space $({\Omega},{\mathcal H},{\mathcal E})$ such that
  ${\mathcal E}[\xi] = -{\mathcal E}[-\xi]$, and ${\widetilde{\mathcal E}}$ be a nonlinear expectation dominated by
  ${\mathcal E}$. Then
 \begin{equation*}
  {\widetilde{\mathcal E}}[\alpha \xi + \eta] = \alpha {\widetilde{\mathcal E}}[\xi] + {\widetilde{\mathcal E}}[\eta], \ \ \text{ for $\forall \alpha \in{\mathbb R}$}.
 \end{equation*}
\end{lemma}

The following lemma  follows easily
with the \eqref{DOM} condition.

  \begin{lemma}\label{lemmono}
    Given a nonlinear expectation space $({\Omega},{\mathcal H},{\widetilde{\mathcal E}})$ and $\xi \in {\mathcal H}$. If $\{{\varphi}_n\}_{n=1}^\infty \subset \mathcal{C}({\mathbb R}^d)$ satisfying ${\varphi}_n\downarrow 0$, then ${\widetilde{\mathcal E}}[{\varphi}_n(\xi)]\downarrow 0$.
  \end{lemma}

\subsection{Martingale problem and the solution}\label{sec:mp}





\begin{definition} [${\widetilde{\mathcal{E}}}$-martingale problem]
  Given the sample space ${\Omega} = C([0,\infty);{\mathbb R}^d)$ with the canonical process $Z$, the {\it ${\widetilde{\mathcal{E}}}$-martingale problem} is to find a time-consistent nonlinear expectation ${\widetilde{\mathcal E}}$ defined on $({\Omega}, {\mathcal H})$ such that,
  for each ${\varphi} \in C_0^\infty({\mathbb R}^d)$,
  \begin{equation*}
  {\varphi}(Z_t) - {\varphi}(Z_0) - \int_0^t {\widetilde{G}}(Z_{\theta},D_z{\varphi}(Z_{\theta}),D_{z}^2{\varphi}(Z_{\theta}))\,d\theta
  \end{equation*}
  is an ${\widetilde{\mathcal{E}}}$-martingale on $[0,\infty)$.
\end{definition}

Now,  we will solve the martingale problem with $Z_t$ in the canonical space ${\Omega} = C_z([0,\infty);{\mathbb R}^{2d})$ being the generalized $G$-Brownian motion in \cite{P10a}. To this end,  consider the following Cauchy problem, which is a special form of PDE \eqref{GE2}.
\begin{equation*}
\begin{cases}
  {\partial}_t u(t,x,y) - {\widetilde{G}}(x,y,D_y u(t,x,y),D_x^2u(t,x,y)) = 0, (t,x,y)\in(0,T]\times{\mathbb R}^d\times{\mathbb R}^d,\\
  u(0,x,y) = \varphi(x,y), (x,y)\in {\mathbb R}^d\times{\mathbb R}^d,
\end{cases}\tag{$\widetilde{\text{P}}$-2}\label{prob2}
\end{equation*}
where $\varphi\in C({\mathbb R}^d\times{\mathbb R}^d)$ with polynomial growth and ${\widetilde{G}}$ satisfies the continuity condition:
\begin{eqnarray}\label{continuity}
|{\widetilde{G}}(z,p,A) - {\widetilde{G}}(\bar{z},p,A)| \le C(1+|A|)(1+|z|^l + |\bar{z}|^l)[(1 + |p|)|z-\bar{z}|]^{\alpha},\quad z\in{\mathbb R}^{2d},
\end{eqnarray}
for some constants $C>0,l\in{\mathbb{N}}$, and ${\alpha} \in(0,1]$.

Clearly the existence and uniqueness of PDE \eqref{GE2} imply the existence and uniqueness of PDE \eqref{prob2}, and both Example~\ref{example1} and Example~\ref{exoftG} satisfy  condition \eqref{continuity}. Note also the $G$ for the $G$-Brownian motion in \cite{P10a} is a special case of the $G$ in the PDE \eqref{PDE}.

\begin{theorem}[Martingale problem] \label{thmmartingale}
   Take the canonical process $(X_t,y_t)$ as the generalized $G$-Brownian motion. Then there exists a time consistent nonlinear expectation ${\widetilde{\mathcal E}}$,
   together with its conditional expectatations $\{{\widetilde{\mathcal E}}_t\}_{t\ge 0}$, defined on the sublinear expectation space $({\Omega}, {\mathcal H})$ such that for $0 \le s \le t \le T$, and ${\varphi}\in C_0^\infty({\mathbb R}^d\times{\mathbb R}^d)$,
  \begin{eqnarray}
  {\widetilde{\mathcal E}}_s[{\varphi}(X_t,y_t)-{\varphi}(X_s,y_s) - \int_s^t {\widetilde{G}}\left(X_{\theta},y_{\theta},D_y{\varphi}(X_{\theta},y_{\theta}), D_{x}^2{\varphi}(X_{\theta},y_{\theta})\right)d\theta] = 0.\label{martingale2}
  \end{eqnarray}
  That is,
    \begin{equation*}
    \begin{aligned}
      {\varphi}(X_t,y_t) - {\varphi}(X_0,y_0) - \int_0^t {\widetilde{G}}\left(X_{\theta},y_{\theta},D_y{\varphi}(X_{\theta},y_{\theta}),D_x^2{\varphi}(X_{\theta},y_{\theta})\right)\,d\theta
    \end{aligned}
    \end{equation*}
    is an ${\widetilde{\mathcal{E}}}$-martingale on $[0,T]$.
\end{theorem}

To prove Theorem \ref{thmmartingale}, it suffices to prove the following Proposition \ref{prop2}.
 Indeed the following
identity can be easily established by taking the It\^o's formula for the generalized $G$-Brownian motion:
  \begin{equation*}
  {\varphi}(X_t,y_t) - {\varphi}(X_s,y_s) = \int_s^t \left\{ \langle D_x{\varphi}(X_{\theta},y_{\theta}),dX_{\theta}\rangle + \langle D_y{\varphi}(X_{\theta},y_{\theta}),dy_{\theta}\rangle +{\frac{1}{2}} \text{tr}[D_x^2{\varphi}(X_{\theta},y_{\theta})\,d{\langle X \rangle}_{\theta}]\right\}.
  \end{equation*}

\begin{proposition}\label{prop2}
  Let $M_0\in{\mathbb R}, \zeta, q \in M_{\mathcal E}^2(0,T;{\mathbb R}^d)$, and $\eta \in M_{\mathcal E}^2(0,T;\mathbb{S}^d)$ be given continuous processes, and let
  \begin{equation*}
  M_t = M_0 + \int_0^t \left(\zeta^T_{\theta}\,dX_{\theta} + q^T_{\theta}\,dy_{\theta} + \text{\em tr}[\eta_{\theta}\,d{\langle X \rangle}_{\theta}]\right) - \int_0^t {\widetilde{G}}(X_{\theta},y_{\theta},q_{\theta},2\eta_{\theta})\,d\theta,\; 0\le t \le T.
  \end{equation*}
  Then $M$ is a ${\widetilde{G}}$-martingale.
\end{proposition}

We will prove Proposition \ref{prop2} in several steps.

\paragraph{Step 1.}

Since $X$ is a symmetric $G$-Brownian motion, by Lemma
\ref{PropI.3.7} and Proposition III.9.1-(iii) of~\cite{P10a}, it is
also a symmetric ${\widetilde{G}}$-martingale. Thus it suffices to prove
that
\begin{equation*}
\bar{M}_{t}=M_{0}+\int_{0}^{t}\left(q^T_{\theta }\,dy_{\theta }+\text{tr}[\eta _{\theta }\,d{%
\langle X\rangle }_{\theta }]\right)-\int_{0}^{t}{\widetilde{G}}(X_{\theta
},y_{\theta },q_{\theta },2\eta _{\theta })\,d\theta ,\quad 0\leq t\leq T
\end{equation*}%
is a ${\widetilde{G}}$-martingale.

To this end, we need the following lemma.

\begin{lemma}
\label{BasicLem} Let $\varphi \in C_{l.Lip}(\mathbb{R}^{2d})$ be given and assume that $f:\mathbb{R}^{2d}\to\mathbb{R}$ satisfies
\[
|f(z)-f(\bar{z})| \le C(1+|z|^l + |\bar{z}|^l)|z-\bar{z}|^\alpha,
\]
for some constants $C>0, l\in\mathbb{N},$ and $\alpha\in(0,1]$. We
have
\begin{equation}
u(t,Z_{t})=\widetilde{\mathcal{E}}_{t}[\varphi
(Z_{T})+\int_{t}^{T}f(Z_{s})ds],  \label{phi-f}
\end{equation}%
where $u\in C([0,T]\times \mathbb{R}^{2d})$ with polynomial growth is the unique viscosity solution of the problem
\begin{equation}
\begin{aligned}
  \partial _{t}u+\widetilde{G}(z,D_{y}u,D_{x}^{2}u)+f(x,y) &=0,\  \ t\in
\lbrack 0,T),\ z=(x,y)\in \mathbb{R}^{2d},\   \\
u(T,x,y) &=\varphi (x,y).
\end{aligned}\label{phi-f-PDE}
\end{equation}
\end{lemma}

\noindent{\textbf{Proof.} }
For a fixed $\bar{t}\in[0,T)$, we set $t_{i}^{n}=i(T-\bar{t})/n$, for $i=0,1,\dots,n,$ and $f_{n}(s,\omega
)=\sum_{i=0}^{n-1}f(\omega(t_{i}^{n}))1_{[t_{i}^{n},t_{i+1}^{n})}(s)$, then
denote
\begin{equation*}
u_{i}^{n}(t,Z_{t};Z_\cdot):=\widetilde{\mathcal{E}}_{t}[\varphi
(Z_{T})+\int_{t}^{T}f_{n}(s,Z_{\cdot })ds],\  \  \ t\in [t_{i}^n,t_{i+1}^n).
\end{equation*}%
According to the definition of the conditional expectation $\widetilde{%
\mathcal{E}}_{t}$, it is not hard to see that $u_{i}^{n}(t,z;\omega)$ solves the following
PDEs parameterized by $\omega$
\begin{align*}
\partial _{t}u_{i}^{n}(t,z;\omega)+\widetilde{G}%
(z,D_{y}u_{i}^n,D_{x}^{2}u_{i}^n)+f(\omega({t_i^n}))& =0,\ t\in [t_{i}^n,t_{i+1}^n),\ z\in
\mathbb{R}^{2d}, \\
\ u_{i}^{n}(t_{i+1}^{n},z;\omega)& =u_{i+1}^{n}(t_{i+1}^{n},z;\omega),
\end{align*}%
for $i=n-1,n-2,\dots ,1,0.$
The terminal condition for $u_k^n$, $k=n-1$,  is at $t_{n}^{n}=T$, $u_{n-1}^{n}(t_{n}^{n},z;\omega)=\varphi (z)$. By the comparison theorem of PDE, backwardly and
successively, we can check that%
\begin{equation*}
u_{i}^{n}(t,z;\omega)-u(t,z)\leq \hat{u}_i^{n}(t,z;\omega%
),\  \  \ t\in [t_{i}^{n},t_{i+1}^{n}),\  \ z\in \mathbb{R}^{2d},
\end{equation*}%
where $\hat{u}_{i}^{n}(t,z;\omega)$, $i=n-1,n-2,\dots,1,0,$ solves the
PDEs
\begin{align*}
\partial _{t}\hat{u}_{i}^{n}(t,z;\omega)+G(D_{z}\hat{u}_i^n,D_{z}^{2}\hat{u}_i^{n})+f(\omega(t_i^n))-f(z)& =0,\  \  \ t\in
[t_{i}^n,t_{i+1}^n), z\in \mathbb{R}^{2d}, \\
\hat{u}_i^{n}(t_{i+1}^{n},z;\omega)& =\hat{u}_{i+1}^{n}(t_{i+1}^{n},z;\omega).
\end{align*}%
Now, since $G(p,A)$ is a sublinear function which does not depend on $z$, we
claim that
\begin{equation*}
|\hat{u}_0^{n}(\bar{t},Z_{\bar{t}};Z_{\cdot})|\le\mathcal{E}_{\bar{t}}[\int_{%
\bar{t}}^{T}|f(Z_{s})-f_{n}(s,Z_{\cdot })|ds]\rightarrow 0 \text{ in $L^1_{\mathcal{E}}(\Omega_T)$ as $n\to\infty$},
\end{equation*}%
from which we have
\begin{equation*}
|u(\bar{t},Z_{\bar{t}})-u_{0}^{n}(\bar{t},Z_{\bar{t}};Z_{\cdot%
})|\leq |\hat{u}_0^{n}(\bar{t},Z_{\bar{t}};Z_{\cdot})|\rightarrow 0\  \  \text{%
as}\  \ n\rightarrow \infty ,
\end{equation*}%
and thus
\begin{eqnarray*}
u_{0}^{n}(\bar{t},Z_{\bar{t}};Z_{\cdot}) &:&=\widetilde{\mathcal{E}}%
_{t}[\varphi (Z_{T})+\int_{t}^{T}f_{n}(s,Z_{\cdot })ds] \\
&\rightarrow &u(\bar{t},Z_{\bar{t}})=\widetilde{\mathcal{E}}%
_{t}[\varphi (Z_{T})+\int_{t}^{T}f(s,Z_{s})ds].
\end{eqnarray*}%

Now let us prove the claim. Set
\[
\Delta = \mathcal{E}_{\bar{t}}[\int_{\bar{t}}^T |f(Z_s)-f_n(s,Z_\cdot)|\,ds],
\]
then
\begin{equation*}
  \Delta \le \mathcal{E}_{\bar{t}}\left[\int_{\bar{t}}^T \sum_{i=0}^{n-1}C(1+|Z_s|^l+|Z_{t_i^n}|^l)|Z_s-Z_{t_i^n}|^\alpha\textbf{1}_{[t_i^n,t_{i+1}^n)}(s)\,ds\right],
\end{equation*}
and
\begin{equation*}
\begin{aligned}
  \mathcal{E}[\Delta] &\le \mathcal{E}\left[\int_{\bar{t}}^T \sum_{i=0}^{n-1}C(1+|Z_s|^l+|Z_{t_i^n}|^l)|Z_s-Z_{t_i^n}|^\alpha\textbf{1}_{[t_i^n,t_{i+1}^n)}(s)\,ds\right]\\
  &\le  C\sum_{i=0}^{n-1}\int_{t_i^n}^{t^n_{i+1}}\mathcal{E}[(1+|Z_s|^l+|Z_{t_i^n}|^l)|Z_s-Z_{t_i^n}|^\alpha]\,ds\\
  &\le  C\sum_{i=0}^{n-1}\int_{t_i^n}^{t^n_{i+1}}\sqrt{\mathcal{E}[(1+|Z_s|^l+|Z_{t_i^n}|^l)^2]\mathcal{E}[|Z_s-Z_{t_i^n}|^{2\alpha}]}\,ds\\
  &\le
  C\sum_{i=0}^{n-1}\int_{t_i^n}^{t^n_{i+1}}(s-t_i^n)^\alpha\,ds\\
  &\le
  C\sum_{i=0}^{n-1}(t_{i+1}^n-t_i^n)^{1+\alpha}\to 0, \text{ as $n\to\infty$}.
\end{aligned}
\end{equation*}
And a slight modification of the approach in Chapter V.3 of \cite{P10a} yields the inequality in the claim.
\endofproof

\paragraph{Step 2:}
Now, we prove that, for each \ fixed $(\eta ,q)\in L_{ip}({\Omega }_{s},{\mathbb{S}}^{d}\times {%
\mathbb{R}}^{d})$, we have the relation, for $0\leq s<t\leq T$,
\begin{equation}\label{2.4.4-00}
\frac{1}{2}\left \langle \eta X_{s},X_{s}\right \rangle +\left \langle
q,y_{s}\right \rangle =\widetilde{\mathcal{E}}_{s}\left[\frac{1}{2}\left \langle
\eta X_{t},X_{t}\right \rangle +\left \langle q,y_{t}\right \rangle
-\int_{s}^{t}\widetilde{G}(X_{\theta },y_{\theta },q,\eta )\,d\theta \right]
\end{equation}%
or
\begin{equation}
\widetilde{\mathcal{E}}_{s}\left[ \frac{1}{2}\text{tr}[\eta (\left \langle
X\right \rangle _{t}-\left \langle X\right \rangle _{s})]
+\left \langle q,y_{t}-y_{s}\right \rangle -\int_{s}^{t}\widetilde{G}(X_{\theta
},y_{\theta },q,\eta )\,d\theta \right]=0.\label{2.4.4-0}
\end{equation}

\noindent{\textbf{Proof of \eqref{2.4.4-0}.}}
We can fix $(\eta ,q)$ as constants. According to  Lemma \ref{BasicLem}, the
right hand side of \eqref{2.4.4-00} equals $u(s,X_{s},y_{s})$, where $u(t,x,y)$ is the viscosity solution of the following PDE:
\begin{align*}
\partial _{t}u+\widetilde{G}(x,y,D_{y}u,D_{x}^{2}u)-\widetilde{G}(x,y,q,\eta)& =0,\ t\in \lbrack
t_{1},t_{2}),\ x,y\in \mathbb{R}^{d}, \\
u(t_{2},x,y)& =\frac{1}{2}\left \langle \eta x,x\right \rangle +\left \langle
q,y\right \rangle .
\end{align*}%
But it is easy to check that $u(t,x,y)\equiv \frac{1}{2}\left \langle \eta
x,x\right \rangle +\left \langle p,y\right \rangle $ is the unique solution of
this PDE, from which we prove the first relation \eqref{2.4.4-00}. For relation~\eqref{2.4.4-0},
we just need to move the terms of the right hand side to the left, inside
the $\widetilde{\mathcal{E}}_{s}$. \ Note  by It\^{o}'s
formula,
\begin{equation*}
\frac{1}{2}\left \langle \eta X_{t},X_{t}\right \rangle -\frac{1}{2}%
\left \langle \eta X_{s},X_{s}\right \rangle =\int_{s}^{t}\langle \eta X_{\theta
},dX_{\theta}\rangle +\frac{1}{2}\text{tr}[\eta (\left \langle X\right \rangle
_{t}-\left \langle X\right \rangle _{s})]
\end{equation*}%
and, for each $\xi \in L_{\mathcal{E}}^{1}(\Omega _{T})$, we have $%
\widetilde{\mathcal{E}}_{s}[\int_{s}^{t}2\langle \eta X_{\theta },dX_{\theta }\rangle+\xi ]=%
\widetilde{\mathcal{E}}_{s}[-\int_{s}^{t}2\langle \eta X_{\theta },dX_{\theta }\rangle+\xi ]=%
\widetilde{\mathcal{E}}_{s}[\xi ]$,
by Lemma~\ref{PropI.3.7}.
\endofproof

\hspace{.2in}

\noindent{\bf Step 3}: Now we are ready to finish proving Proposition \ref{prop2}.

By the domination inequality, the equality~\eqref{2.4.4-0} can be extended to the case
$(\varsigma,\eta,q)\in L_{\mathcal E}^{2}({\Omega}_{s},\mathbb{R}\times
{\mathbb{S}}^{d}\times{\mathbb{R}}^{d})$. Now for step processes:
\begin{align*}
\eta_{t}^{K}=\sum_{j=0}^{K-1}\eta_{j}\mathbf{1}_{[s_{j},s_{j+1})}%
(t),\ \ q_{t}^{K} &  =\sum_{j=0}^{K-1}q_{j}\mathbf{1}_{[s_{j},s_{j+1}%
)}(t),\ s=s_{1}<\cdots<s_{K}=T\text{,}\ \\
(\eta_{j},q_{j}) &  \in L_{\mathcal E}^{2}({\Omega}_{s_{j}},{\mathbb{S}}^{d}%
\times{\mathbb{R}}^{d}),\ \varsigma\in L_{\mathcal E}^{2}({\Omega}_{s},\mathbb{R}).
\end{align*}
We can repeat the equality \eqref{2.4.4-0} to prove
\[
\widetilde{\mathcal{E}}_{s}\left[\varsigma+\int_{s}^{t}\text{tr}[\eta_{r}^{K}d\langle
X\rangle_{r}]+\int_{s}^{t}\langle q_{r}^{K},dy_{r}\rangle -\int_{s}^{t}\widetilde{G}%
(r,{\omega},q_{r}^{K},2\eta_{r}^{K})dr\right]=\varsigma,\ \ s\leq t\leq T.
\]
From the domination of ${\widetilde{\mathcal E}}$ by $\mathcal{E}$, we then can prove, for $(\eta_{\cdot},q_{\cdot})\in
M_{\mathcal E}^{2}(0,T;{\mathbb{S}}^{d}\times{\mathbb{R}}^{d})$, and
\[
\widetilde{\mathcal{E}}_{s}\left[\int_{s}^{t}\text{tr}[\eta_{r}d\langle X\rangle_{r}]+\int%
_{s}^{t}\langle q_{r},dy_{r}\rangle -\int_{s}^{t}\widetilde{G}({\omega}(s),q_{r},2\eta_{r})dr\right]=0,
\]
from which the proof is complete.\endofproof

\hspace{.2in}


  As a corollary of Theorem~\ref{thmmartingale}, we have the following result.

  \begin{corollary}\label{corofmp}
    The martingale Theorem~\ref{thmmartingale} holds when ${\varphi}$ is a polynomial.
  \end{corollary}

  \noindent{\bf Proof.} For each given polynomial ${\varphi}$, one can find a sequence of functions ${\varphi}_n \in C_0^\infty({\mathbb R}^d\times{\mathbb R}^d)$ such that $|{\varphi}_n -{\varphi}| \downarrow 0, |D_y{\varphi}_n -D_y{\varphi}|\downarrow 0,$ and $|D_x^2{\varphi}_n -D_x^2{\varphi}|\downarrow 0$, then we have
  \begin{align*}
    &\,\bigg|{\widetilde{\mathcal E}}_s[\varphi(Z_t)-\varphi(Z_s) - \int_s^t {\widetilde{G}}(Z_{\theta},D_y\varphi(Z_{\theta}),D_x^2\varphi(Z_{\theta}))\,d\theta]\bigg|\\
    \le&\,\bigg|{\widetilde{\mathcal E}}_s[\varphi(Z_t)-\varphi(Z_s) - \int_s^t {\widetilde{G}}(Z_{\theta},D_y{\varphi}(Z_{\theta}),D_x^2{\varphi}(Z_{\theta}))\,d\theta] \\
    &\,- {\widetilde{\mathcal E}}_s[\varphi_n(Z_t)-\varphi_n(Z_s) - \int_s^t {\widetilde{G}}(Z_{\theta},D_y{\varphi}_n(Z_{\theta}),D_x^2{\varphi}_n(Z_{\theta}))\,d\theta]\bigg|\\
    \le&\,\bigg|{\mathcal E}_s[(\varphi(Z_t)-\varphi_n(Z_t))- (\varphi(Z_s)-\varphi_n(Z_s))\\
     &\,- \int_s^t [{\widetilde{G}}(Z_{\theta},D_y{\varphi}(Z_{\theta}),D_x^2{\varphi}(Z_{\theta}))- {\widetilde{G}}(Z_{\theta},D_y{\varphi}_n(Z_{\theta}),D_x^2{\varphi}_n(Z_{\theta}))]\,d\theta]\bigg|\\
    \le&\,{\mathcal E}_s[|\varphi(Z_t)-\varphi_n(Z_t)|]+{\mathcal E}_s[|\varphi(Z_s)-\varphi_n(Z_s)|]\\
    &\,+\int_s^t{\mathcal E}_s[|{\widetilde{G}}(Z_{\theta},D_y{\varphi}(Z_{\theta}),D_x^2{\varphi}(Z_{\theta}))-{\widetilde{G}}(Z_{\theta},D_y{\varphi}_n(Z_{\theta}),D_x^2{\varphi}_n(Z_{\theta}))|]\,d\theta\\
    \le&\,{\mathcal E}_s[|\varphi(Z_t)-\varphi_n(Z_t)|]+{\mathcal E}_s[|\varphi(Z_s)-\varphi_n(Z_s)|]\\
    &\,+ L\int_s^t{\mathcal E}_s[|D_y{\varphi}(Z_{\theta})-D_y{\varphi}_n(Z_{\theta})|+|D_x^2{\varphi}(Z_{\theta})-D_x^2{\varphi}_n(Z_{\theta})|]\,d\theta.
  \end{align*}
  Letting $n\to\infty$, according to Lemma \ref{lemmono}, we see
  \begin{equation*}
  {\widetilde{\mathcal E}}_s[{\varphi}(Z_t)-{\varphi}(Z_s) - \int_s^t {\widetilde{G}}(Z_{\theta},D_y{\varphi}(Z_{\theta}),D_x^2{\varphi}(Z_{\theta}))\,d\theta] = 0.
  \end{equation*}
  \endofproof



\section{Weak solution of $G$-SDE}\label{sec:weak}
In this section, we will develop a notion of weak solution of general $G$-SDE, and show the existence of such weak solutions in comparison with the strong solutions within the existing $G$-framework. To this end, the $G$ in this section is from \cite{P10a} which is a special case of the $G$ (with the form~\eqref{expofG1}) in  class $\mathfrak{D}$.  We will rely on the analysis and results for the martingale problem in the previous section.

\subsection{Weak solution of SDE in  $G$-framework}

Consider a $d$-dimensional $G$-SDE
\begin{eqnarray}\label{GSDE}
\begin{cases}
  dz_t = b(z_t)dt + r(z_t)d{\langle B \rangle}_t + {\sigma}(z_t) dB_t,\quad 0\le t \le T,\\
  z_0 = z, \quad z\in{\mathbb R}^d,
\end{cases}
\end{eqnarray}
where $b = (b^i)_{1\le i \le d},r = (r_{jk}^i)_{1\le i,j,k \le d}$, and ${\sigma} = ({\sigma}_{ij})_{1\le i,j \le d}$, and $b^i,r_{jk}^i,$ and ${\sigma}_{ij}$ are continuous functions on ${\mathbb R}^d$, and ${\sigma}\ge{\sigma}_0 I$ for some constant ${\sigma}_0>0$. Here $B$ is a $d$-dimensional generalized ${\mathcal E}$-Brownian motion introduced in~\cite{P10a} (for additional background material, please consult Appendix).


\begin{definition}[Weak solution]\label{DefofSol}
  A {\it weak solution} of $G$-SDE~\eqref{GSDE} is a triple $(({\Omega},{\mathcal H},{\widetilde{\mathcal E}}),z,B)$, where
  \begin{enumerate}[i)]
    \item $({\Omega},{\mathcal H},{\widetilde{\mathcal E}})$ is a nonlinear expectation space,
    \item $z$ is a $d$-dimensional continuous process on the nonlinear expectation space $({\Omega},{\mathcal H},{\widetilde{\mathcal E}})$, and $B$ is a $d$-dimensional ${\widetilde{\mathcal E}}$-Brownian motion in the sense of \cite{P10a},
    \item the identity
    \begin{eqnarray}
    z_t = z_0 + \int_0^t b(z_{\theta})\,d\theta + \int_0^t r(z_{\theta}) \,d{\langle B \rangle}_{\theta} +  \int_0^t {\sigma}(z_{\theta})\,d B_{\theta}
    \end{eqnarray}
    holds in the nonlinear expectation space.
  \end{enumerate}
\end{definition}

\subsection{Existence of weak solutions for $G$-SDE}

We will establish the existence result for a random process $(z_{t})_{t\geq0}$
which is a weak solution of a $G$-SDE. For comparison with the existing
$G$-framework, we assume that a random process $z_{t}$ in a nonlinear
expectation space can be decomposed into two parts $X_{t}$ and $y_{t}$, with
$X_{t}$ being a symmetric martingale and $y_{t}$ with finite variation. This
is a generalization of the generalized $G$-Brownian motion $(X_{t},y_{t})$. For more
about this martingale representation in a nonlinear expectation space, see
for instance \cite{PSZ12}. More specifically, we consider the following
1-dimensional $G$-SDE:
\begin{equation}%
\begin{cases}
dz_{t}^{i}=b^{i}(z_{t})\,dt+r_{jk}^{i}(z_{t})\,d{\langle B}{\rangle}^{jk}%
_{t}+{\sigma}_{ij}(z_{t})\,dB_{t}^{j},\quad0\leq t\leq T,\\
z_{0}=z\in{\mathbb{R}}^{d},
\end{cases}
\label{SGSDE}%
\end{equation}
where $b:{\mathbb{R}}^{d}\to{\mathbb{R}}^{d}$, $r:\mathbb{R}^{d}\to
L(\mathbb{R}^{d\times d};\mathbb{R}^d),$ and $\sigma:\mathbb{R}^{d}%
\to L(\mathbb{R}^{d};\mathbb{R}^d)$ are bounded and continuous functions such that the inversed matrix ${\sigma}^{-1}(z)$ is also bounded, and they satisfy the H\"older continuity condition
\[
|b(z)-b(\bar{z})|+|r(z)-r(\bar{z})|+|{\sigma}(z)-{\sigma}(\bar{z})|\leq
L_{0}|z-\bar{z}|^{\alpha}.
\]
We also assume that there exists a sublinear monotone function $\bar{G}:\mathbb{S}^{d}%
\mapsto \mathbb{R}$, satisfying
\[
\bar{G}(A)\geq \lambda \text{tr}[A],\  \  \forall A\in \mathbb{S}^{d}%
\  \ (\lambda>0),
\]
define
\[
{\widetilde{G}}(x,y,p,A)=\bar{G}([2r_{ij}^{k}(x+y)p_{k}+{\sigma}_{i^{\prime}i%
}(x+y)\sigma_{j^{\prime}j}(x+y)A^{i'j'}]_{i,j=1}^{d})+b^{i}(x+y)p_{i},
\]
here we use Einstein convention, namely the repeated indices $i,j$ implies
taking sum from $1$ to $d$. One can see that $\widetilde{G}$ satisfies the continuity condition~\eqref{continuity}.

\begin{lemma}
\label{N-p-eta}For the case
\[
{\widetilde{G}}(x,y,p,A)=\bar{G}([2r_{ij}^{k}(x+y)p_{k}+{\sigma}_{i^{\prime}i%
}(x+y)\sigma_{j^{\prime}j}(x+y)A^{i'j'}]_{i,j=1}^{d})+b^{i}(x+y)p_{i},%
\]
we denote $z_{t}=X_{t}+y_{t}$, where $(X_{\cdot},y_{\cdot})=\omega({\cdot})%
\in \Omega$ is the canonical process. Then, for each $p\in M_{\bar{\mathcal{E}}}^{2}(0,T;\mathbb{R}^d)$ and
$\eta \in M_{\bar{\mathcal{E}}}^{1}(0,T;
\mathbb{S}^d)$, the process
\begin{equation}
N_{t}^{p,\eta}=\int_{0}^{t}p_{s}^T \,dz_{s}+\int_{0}^{t}\frac{1}{2}\text{\em tr}[\eta
_{s}\,d\left \langle z\right \rangle _{s}]-\int_{0}^{t}[\bar{G}(2r(z_{s}%
)p_{s}+\sigma^{T}(z_s)\eta_{s}\sigma(z_{s}))+p_{s}^Tb(z_{s}%
) ]\,ds \label{G-Z-martingale}%
\end{equation}
is a martingale under $\widetilde{\mathcal{E}}$.
\end{lemma}

\noindent{\bf Proof.}
Since, for each $\xi,p$, and $\eta$
\[
\int_{0}^{t}\zeta_{s}^T\,dX_{s}+\int_{0}^{t}p_{s}^T\,dy_{s}+\int_{0}^{t}\frac{1}%
{2}\text{tr}[\eta_{s}\,d\left \langle X\right \rangle _{s}]-\int_{0}^{t}\widetilde{G}(X_{s}%
,y_{s},p_{s},\eta_{s})ds
\]
is a martingale under $\widetilde{\mathcal{E}}$ and $\left \langle
X\right \rangle _{s}\equiv \left \langle z\right \rangle _{s}$, by taking
$\zeta_{s}\equiv p_{s}$ we obtain that for each $p\in M_{\bar{\mathcal{E}}}^{2}(0,T;\mathbb{R}^d)$ and
$\eta \in M_{\bar{\mathcal{E}}}^{1}(0,T;\mathbb{S}^d)$, $N_{t}^{p,\eta}$ is an $\widetilde{\mathcal{E}}$-martingale.
\endofproof

\begin{theorem}
\label{Existence} Under the nonlinear expectation ${\widetilde{\mathcal{E}}}$
derived from the PDE
\[
\partial_{t}u(t,x,y)-{\widetilde{G}}(x,y,D_{y}u,D_{x}^{2}u)=0,\;(t,x,y)\in
(0,T]\times{\mathbb{R}}^{2d},
\]
together with the canonical space ${\Omega}=C_{x,y}([0,T];{\mathbb{R}}^{2d})$,
the linear space of d-dimensional random variables ${\mathcal{H}}$, the process $z_{\cdot}=X_{\cdot}+y_{\cdot}$ for $(X_{\cdot},y_{\cdot
})\in{\Omega}$, and
\begin{align*}
B_{t}^{i}   =\int_{0}^{t}\sigma_{ij}^{-1}(z_{s})dz_{s}^{j}-\int_{0}%
^{t}\sigma_{ij}^{-1}(z_{s})b^{j}(z_{s})ds -& \int_{0}^{t}\sigma_{ik}^{-1}%
(z_{s})r_{jl}^{k}(z_{s})\,d\left \langle B\right \rangle _{s}^{jl},\  \\
&  \qquad i  =1,\dots,d,\  \  \ 0\leq t\leq T,
\end{align*}
is a weak solution of the $G$-SDE~\eqref{SGSDE}.
\end{theorem}

Several steps are needed for proving the theorem.

\noindent{\bf Step 1.}
For the canonical process $(z_{s})_{s\geq0}$, we construct the following
It\^{o} process:%
\begin{equation}
B_{t}^{i}=\int_{0}^{t}\sigma_{ij}^{-1}(z_{s})dz_{s}^{j}-\int_{0}^{t}%
\sigma_{ij}^{-1}(z_{s})b^{j}(z_{s})ds-\int_{0}^{t}\sigma_{ik}^{-1}%
(z_{s})r_{jl}^{k}(z_{s})d\left \langle B\right \rangle _{s}^{jl}.
\label{B1}%
\end{equation}
By Proposition \ref{Quadra-xi}, the quadratic variation process of this $G$-It\^{o}
process $B$ is given by
\begin{equation}
\left \langle B\right \rangle _{t}^{ii^{\prime}}=\int_{0}^{t}\sigma_{ij}%
^{-1}(z_s)\sigma_{i^{\prime}j^{\prime}}^{-1}(z_{s})d\left \langle z\right \rangle_s
^{jj^{\prime}}. \label{Quad-B}%
\end{equation}
Thus it is clear that
\[
dz_{s}=b(z_{s})ds+\sigma(z_{s})dB_{s}+r(z_{s})d\left \langle
B\right \rangle _{s},\  \  \ z_0= z (= x+y).
\]
It remains to prove that, under $\widetilde{\mathcal{E}}$, $B$ is a $d$-dimensional $\bar{G}$-Brownian
motion. From (\ref{Quad-B}), we rewrite (\ref{B1}) as%
\begin{align}
B_{t}^{l}  &  =\int_{0}^{t}\sigma_{lj}^{-1}(z_{s})dz_{s}^{j}-\int_{0}%
^{t}\sigma_{lj}^{-1}(z_{s})b^{j}(z_{s})ds\label{B2}\\
&  - \sum_{i,i'}\int_{0}^{t}\sigma_{lk}^{-1}(z_{s})r_{ii^{\prime}}^{k}(z_{s}%
)\sigma_{ij}^{-1}(z_{s})\sigma_{i^{\prime}j^{\prime}}^{-1}(z_{s})d\left \langle
z\right \rangle_s ^{jj^{\prime}}.\nonumber
\end{align}
We need to prove that the $G$-It\^{o} process defined by
\begin{equation}
M_{t}:=\sum_{k}\int_{0}^{t}\zeta_{s}^{k}dB_{s}^{k}+\frac{1}{2}\sum_{i,i'}\int_{0}^{t}\eta
_{s}^{ii^{\prime}}d\left \langle B\right \rangle _{s}^{ii^{\prime}}-\int_{0}%
^{t}\bar{G}(\eta_{s})ds \label{M}%
\end{equation}
is an $\widetilde{\mathcal{E}}$-martingale. Indeed, we have
\begin{align*}
M_{t}  &  =\int_{0}^{t}\zeta_{s}^{l}\sigma_{lj}^{-1}(z_{s})dz_{s}^{j}\\
&  +\int_{0}^{t}\left[\frac{1}{2}\eta_{s}^{ii^{\prime}}\sigma_{ij}^{-1}(z_{s})%
\sigma_{i^{\prime}j^{\prime}}^{-1}(z_{s})-\sum_{i,i'}\zeta_{s}^{l}\sigma_{lk}^{-1}%
(z_{s})r_{ii^{\prime}}^{k}(z_{s})\sigma_{ij}^{-1}(z_{s})\sigma_{i^{\prime
}j^{\prime}}^{-1}(z_{s})\right]d\left \langle z\right \rangle_s ^{jj^{\prime}}\\
&  -\int_{0}^{t}[\bar{G}(\eta_{s})+\zeta_{s}^{l}\sigma_{lj}^{-1}(z_{s}%
)b^{j}(z_{s})]ds=N_{t}^{\bar{p},\bar{\eta}},
\end{align*}
where we set
\begin{align*}
(\bar{p}_{s})_{j}  &  =\zeta_{s}^{l}\sigma_{lj}^{-1}(z_{s}), \text{ and}\\
(\bar{\eta}_{s})_{jj^{\prime}}  &  =\sigma_{ij}^{-1}(z_{s})\eta_{s}^{ii^{\prime
}}(z_{s})\sigma_{i^{\prime}j^{\prime}}^{-1}(z_{s})-2\sum_{i,i'}\zeta_{s}^{l}\sigma_{lk}%
^{-1}(z_{s})r_{ii^{\prime}}^{k}(z_{s})\sigma_{ij}^{-1}(z_{s})\sigma_{i^{\prime
}j^{\prime}}^{-1}(z_{s}).
\end{align*}
This, with Lemma \ref{N-p-eta}, shows that $\{M_{t}\}_{0\le t\le T}$ defined in
(\ref{M}) is an $\widetilde{\mathcal{E}}$-martingale.

\hspace{.1in}

\noindent{\bf Step 2.}
The proof can be completed by applying the following proposition.

\begin{proposition}
\label{Prop-G-mart}Let $(B_{t})_{t\geq0}$ be a $d$-dimensional $G$-It\^{o}
 process defined on sublinear expectation space $(\Omega,L_{\mathcal{E}%
}^{1}(\Omega),\mathcal{E})$ and let $\bar{G}$ be dominated by $G$ such that for each
bounded $\zeta \in M_{\mathcal{E}}^{2}(0,T;\mathbb{R}^{d})$ and $\eta \in M_{\mathcal{E}}%
^{2}(0,T;\mathbb{S}^{d})$,
\[
\int_{0}^{t}\zeta^T_{s}\,dB_{s}+\frac{1}{2}\int_{0}^{t}\text{\em tr}[\eta_{s}d\left \langle
B\right \rangle _{s}]-\int_{0}^{t}\bar{G}(\eta_{s})\,ds,
\]
is an $\widetilde{\mathcal{E}}$-martingale, where $\bar{G}$ given as before.
Then $B$ is a $\bar{G}$-Brownian motion under $\widetilde{\mathcal{E}}$.
\end{proposition}

\begin{corollary}
We assume the same condition for $B_{t}$ and $\bar{G}$ as given in Proposition
\ref{Prop-G-mart}. If for each $\varphi \in C_{b}^{2}(\mathbb{R}^{d})$,%
\[
\int_{0}^{t}\langle D_{x}\varphi(B_{s}),dB_{s}\rangle+\frac{1}{2}\int_{0}^{t}%
\text{\em tr}[D_{x}^{2}\varphi(B_{s})\,d\left \langle B\right \rangle _{s}]-\int
_{0}^{t}\bar{G}(D_{x}^{2}\varphi(B_{s}))ds,
\]
is an $\widetilde{\mathcal{E}}$ martingale, then $B$ is also a $\bar{G}%
$-Brownian motion under $\widetilde{\mathcal{E}}$.
\end{corollary}


\begin{corollary}
We set $d=1$ and assume the same condition for $B$ as in Proposition
\ref{Prop-G-mart}. If there there two constants $\overline{\sigma}%
>\underline{\sigma}>0$, such that
\begin{enumerate}[(1).]
\item $B$ is a symmetric ${\widetilde{\mathcal{E}}}$-martingale,

\item the process $\{B_{t}^{2}-\overline{\sigma}^{2}t\}$ is an ${\widetilde
{\mathcal{E}}}$-martingale,

\item the process $\{ \underline{\sigma}^{2}t-B_{t}^{2}\}$ is an
${\widetilde{\mathcal{E}}}$-martingale.
\end{enumerate}

Then $B$ is a $\bar{G}$-Brownian motion under ${\widetilde{\mathcal{E}}}$,
where $\bar{G}$ is a sublinear monotone function of the form
\[
\bar{G}(a):=\frac{1}{2}(\overline{\sigma}^{2}a^{+}-\underline{\sigma}^{2}%
a^{-}),\  \  \  \ a\in \mathbb{R}\text{.}%
\]
\end{corollary}

\noindent{\bf Proof.}
For simple processes $\zeta$,
$\eta \in M_{\mathcal{E}}^{2,0}(0,T)$ of the form
\[
\zeta_{t}=\sum_{i=0}^{n-1}\zeta^{i}1_{[t_{i},t_{i+1})}(t),\  \  \eta_{t}%
=\sum_{i=0}^{n-1}\eta^{i}1_{[t_{i},t_{i+1})}(t),\  \  \  \zeta^{i},\eta^{i}\in
L_{ip}(\Omega_{t_{i}}),\  \
\]
we can check that the process defined by $M_{t}:=\int_{0}^{t}\zeta_{s}%
dB_{s}+\frac{1}{2}\int_{0}^{t}\eta_{s}d\left \langle B\right \rangle _{s}%
-\int_{0}^{t}\bar{G}(\eta_{s})ds$ is an ${\widetilde{\mathcal{E}}}$-martingale. It is
also easy to extend this property to the case of bounded $\zeta,\eta \in
M_{\mathcal{E}}^{2}(0,T)$. Thus Proposition \ref{Prop-G-mart} can be applied.
\endofproof

\hspace{.1in}

\noindent{\bf Proof of Proposition \ref{Prop-G-mart}. }
For each $\varphi(x,\bar{x})\in
C_{b}^{3}(\mathbb{R}^{2d})$, we solve the following PDE, parameterized by
$\bar{x}\in \mathbb{R}^{d}$,
\begin{equation}
\partial_{t}u_{\varepsilon}(t,x;\bar{x})+\bar{G}(D_{x}^{2}u_{\varepsilon
}(t,x;\bar{x}))=0,\label{PDEe}%
\end{equation}
defined on $t\in \lbrack0,T+\varepsilon)\times \mathbb{R}^{d}$ with terminal
condition $u_{\varepsilon}(T+\varepsilon,x;\bar{x})=\varphi(x,\bar{x})$. Since
$\bar{G}$ is convex and $\bar{G}(A)\geq \lambda$tr$[A]$, by Krylov \cite{K87}, 
the internal regularity
\[
\left \Vert u_\varepsilon\right \Vert _{C^{1+\alpha/2,2+\alpha}([0,T]\times \mathbb{R}^{d}%
)}<\infty.
\]
We then can apply $G$-It\^{o}'s formula to get
\begin{align*}
M_{\bar{t}}^{\varepsilon,t} &  :=u_{\varepsilon}(\bar{t},B_{\bar{t}}^{t};B_{t})-u_{\varepsilon
}(t,0;B_{t})\\
&  =\int_{t}^{\bar{t}}\partial_{t}u_{\varepsilon}(s,B_{s}^{t};B_{t})\,ds +\int_{t}^{\bar{t}}\langle D_{x}u_{\varepsilon}(s,B_{s}^{t};B_{t}),dB_s\rangle +\frac{1}%
{2}\int_{t}^{\bar{t}}\text{tr}[D_{x}^{2}u_{\varepsilon}(s,B_{s}^{t};B_{t})\,d\left \langle
B\right \rangle _{s}]\\
&  =\int_{t}^{\bar{t}}\langle D_{x}u_{\varepsilon}(s,B_{s}^{t};B_{t}),dB_{s}\rangle+\frac
{1}{2}\int_{t}^{\bar{t}}\text{tr}[D_{x}^{2}u_{\varepsilon}(s,B_{s}^{t};B_{t})\,d\left \langle
B\right \rangle _{s}]\\
&  -\int_{t}^{\bar{t}}\bar{G}(D_{x}^{2}u_{\varepsilon}(s,B_{s}^{t};B_{t}))ds,
\end{align*}
where $B^t_s = B_s - B_t, t\le s \le T$.

But, as a condition of Proposition~\ref{Prop-G-mart}, $M^{\varepsilon,t}$ is an ${\widetilde
{\mathcal{E}}}$-martingale. It then follows that
\[
u_{\varepsilon}(t,0;B_t)=\widetilde{\mathcal{E}}_{t}[u_{\varepsilon}(T,B_{T}%
-B_{t};B_{t})].\  \
\]
Let $u$ be the viscosity solution of the same
PDE (\ref{PDEe}) defined on $[0,T)\times \mathbb{R}^{d}$ with terminal value
$u(T,x;\bar{x})=\varphi(x,\bar{x})$. By using the stability of viscosity solution (Lemma II.6.2 of~\cite{FS92}) and the internal regularity of $u$, letting $\varepsilon \rightarrow0$ in the above identity, we obtain $u(t,0;B_t)=\widetilde{\mathcal{E}}_t[u(T,B_{T}-B_{t};B_{t})]=\widetilde{\mathcal{E}}_t%
[\varphi(B_{T}-B_{t},B_{t})]$. It follows that
\[
\widetilde{\mathcal{E}}_{t}[\varphi(B_{T}-B_{t},B_{t})]=\widetilde
{\mathcal{E}}_{t}[\bar{\mathcal{E}}[\varphi(\sqrt{T-t}\xi,\bar{x}%
)]_{\bar{x}=B_{t}}],
\]
where $\xi$ is a $\bar{G}$-normal distributed random variable. It follows that
$B_{T}-B_{t}\overset{d}{=}\sqrt{T-t}\xi$ and $B_{T}-B_{t}$ is independent of
$B_{t}$. In fact, we can applying the above method to the case $\varphi
=\varphi(B_{T}-B_{t},B_{t_{1}},\cdots,B_{t_{N}})$, for $t_{1}\leq \cdots \leq
t_{N}\leq t$, to prove that, for $\varphi(x_{1},\cdots,x_{N},x)\in C_{b}%
^{3}(\mathbb{R}^{d\times(N+1)})$, we have
\begin{align*}
\widetilde{\mathcal{E}}[\varphi(B_{t_{1}},\cdots,B_{t_{N}},B_{T}-B_{t})]  &
=\widetilde{\mathcal{E}}[\widetilde{\mathcal{E}}[\varphi(x_{1},\dots,
x_{N},B_{T}-B_{t})]_{x_{1}=B_{t},\cdots,x_{N}=B_{t_{N}}}]\\
& =\widetilde{\mathcal{E}}[\bar{\mathcal{E}}[\varphi(x_{1},\dots,
x_{N},\sqrt{T-t}\xi)]_{x_{1}=B_{t},\cdots,x_{N}=B_{t_{N}}}].
\end{align*}
This implies that $B_{T}-B_{t}$ is also independent of $B_{t_{1}}%
,\cdots,B_{t_{N}}$. It then follows that $(B_{t})_{t\geq0}$ is a $\bar{G}$-Brownian
motion. The proof is complete.
\endofproof

\begin{remark}
The method to establish the existence of weak solutions of  SDE (\ref{GSDE})  is by and  large a generalization of the classical  Girsanov transformation for change of measures.
However, the Girsanov transformation is limited to the transform of two measures that are absolutely continuous, and even a small change of the diffusion coefficient may cause the singularity between two measures. In this regard, our method is new and the key is to have a sublinear expectation of $\mathcal{E}$ that dominates a class of probability measures singular from each other.

\end{remark}


\section{Appendix}

\subsection{Related stochastic calculus under nonlinear expectations}

We first recall some notions under $G$-framework mainly from \cite{P10a}. We then develop some new results under general sublinear expectations.

\subsubsection{Review: It\^{o}'s integral with $G$-Brownian motion in $L_{\mathcal{E}}%
^{2}(\Omega)$}

We briefly present  some useful results of stochastic calculus under
$G$-expectation. Recall that, since $\mathcal{E}$ is a sublinear expectation
defined on $(\Omega,L_{ip}(\Omega_{T}))$, thus, for each $p\geq1$, $T>0$, we
can define a Banach norm
\[
\left \Vert \xi \right \Vert _{L_{\mathcal{E}}^{p}}=\left(\mathcal{E}[|\xi
|^{p}]\right)  ^{1/p},\  \  \  \xi \in L_{ip}(\Omega_{T}).
\]
The completion of $L_{ip}(\Omega_{T})$ under this norm is denoted by
$L_{\mathcal{E}}^{p}(\Omega_{T})$. Both $\mathcal{E}$-expectation and $\widetilde
{\mathcal{E}}$-expectation, as well as their conditional expectations
$\mathcal{E}_{t}$, $\widetilde{\mathcal{E}}_{t}$ are extended in $L_{\mathcal{E}}%
^{1}(\Omega_{T}), T\geq0$ and the properties obtained in Proposition
\ref{PropofExp} still hold true for $L_{\mathcal{E}}^{p}(\Omega_{T})$ in the place of
$L_{ip}(\Omega_{T})$. Moreover, it is proved in \cite{DHP11} that there exists a
weakly compact subset $\mathcal{P}_{G}$ of probability measures on the Borel
measurable space $(\Omega,\mathcal{B}{\normalsize (\Omega))}$ such that
\[
\mathcal{E}[\xi]=\sup_{P\in \mathcal{P}_{G}}\int_{\Omega}\xi(\omega
)dP,\  \  \xi \in L_{\mathcal{E}}^{1}(\Omega_{T}),
\]
and, in fact $L_{\mathcal{E}}^{p}(\Omega_{T})$ belongs to the space of $\mathcal{B}%
{\normalsize (\Omega)}$-measurable functions
\begin{equation}
\sup_{P\in \mathcal{P}_{G}}\int_{\Omega}|\xi(\omega)|^{p}dP<\infty.\label{Lp}%
\end{equation}
The usual language of $P$-almost surely is replaced by $c_{G}$-quasi surely
with
\[
c_{G}(A):=\sup_{P\in \mathcal{P}_{G}}P(A),\  \ A\in \mathcal{B}%
{\normalsize (\Omega).}%
\]
In fact $\xi \in L_{\mathcal{E}}^{p}(\Omega)$ iff $\xi \in L^{0}(\Omega)$ has $c_{G}%
$-quasi continuous modification such that (\ref{Lp}) and
\[
\lim_{N\rightarrow \infty}\sup_{P\in \mathcal{P}_{G}}\int_{\Omega}|\xi
(\omega)|^{p}1_{\{|\xi|>N\}}dP=0.\text{ }%
\]

We give the definition of $\widetilde{G}$-Brownian motion here.

\begin{definition}[$\widetilde{G}$-Brownian Motion]
A $d$-dimensional process $(B_{t})_{t\geq0}$ defined on a sublinear
expectation space $(\Omega,\mathcal{H},\mathcal{E})$ is called a {\it $\widetilde{G}$-Brownian
motion} under a given nonlinear expectation $\widetilde{\mathcal{E}}${
dominated by }$\mathcal{E}$, if the following conditions are satisfied:
\begin{enumerate}[(i).]
  \item $B_{0}(\omega)${$=0$.
  \item For each
$t,s\geq0$, }$B${{$_{t+s}-B_{t}$ and $B_{s}$ are identically distributed and
}}$B${{{$_{t+s}-B_{t}$ is }independent from $(B_{t_{1}},B_{t_{2}}%
,\cdots,B_{t_{n}})$, for each $n\in \mathbb{N}$ and $0\leq t_{1}\leq \cdots \leq
t_{n}\leq t$.
\item $\lim_{t\downarrow0}\mathcal{E}%
[|B_{t}|^{3}]t^{-1}=0$. }}
\end{enumerate}
$B$ is called a {\it symmetric} Brownian motion if
$\widetilde{\mathcal{E}}[B_{t}]=\widetilde{\mathcal{E}}[-B_{t}]=0$.
In the finite dimensional case nonlinear distribution of $B$ is fully determined
by the function: $\widetilde{G}(A)=\frac{1}{2}\widetilde{\mathcal{E}}[\langle
AB_{1},B_{1}\rangle]$, defined on $\mathbb{S}^{d}$. We often call it a
$\widetilde{G}$-Brownian motion.
\end{definition}

In fact in the main text of the paper we have introduced $2d$-dimensional stochastic
process $(X_{t}-X_{0},y_{t}-y_{0})$. They are Brownian motion under
$\mathcal{E}$, but in general not under $\widetilde{\mathcal{E}}$. Futhermore,
under $\mathcal{E}$, $X_{t}$ is a symmetric $G_{0}$-Brownian motion with
$G_{0}(A)=G(0,A)$, while $y_{t}$ is not symmetric. Indeed, for
each $\varphi \in C_{b}(\mathbb{R}^d)$
\[
\mathcal{E}[\varphi(y_{t})]=\max_{v\in \Gamma_{1}}\varphi(vt),
\]
where $\Gamma_{1}$ is a convex subset of $\mathbb{R}^{d}$ such that
$\max_{v\in \Gamma_{1}}\left \langle p,v\right \rangle =g_{0}(p):=G(p,0)$, for
all $p\in \mathbb{R}^{d}$. A typical situation of such kind of Brownian motion
is the quadratic variation process $\left \langle X\right \rangle _{t}$ of the
above symmetric Brownian motion.


\begin{definition}
\label{d-Def-4}For $T\in [0,\infty),$ a partition $\pi_{T}$ of $[0,T]$ is a
finite ordered subset $\pi_{T}=\{t_{0},t_{1},\dots,t_{N}\}$ such that
$0=t_{0}<t_{1}<\cdots<t_{N}=T$,
\[
\mu(\pi_{T}):=\max \{|t_{i+1}-t_{i}|; i=0,1,\dots,N-1\} \text{.}%
\]
We use $\pi_{T}^{N}=\{t_{0}^{N},t_{1}^{N},\dots,t_{N}^{N}\}$ to denote a
sequence of partitions of $[0,T]$ such that $\displaystyle\lim_{N\rightarrow \infty}\mu
(\pi_{T}^{N})=0$.
\end{definition}

Let $p\geq1$ be fixed. We consider the following type of simple processes: for
a given partition $\pi_{T}=\{t_{0},\dots,t_{N}\}$ of $[0,T]$ we set%
\[
\eta_{t}(\omega)=\sum_{k=0}^{N-1}\xi_{k}(\omega)\textbf{1}_{[t_{k},t_{k+1}%
)}(t),
\]
where $\xi_{k}\in L_{\mathcal{E}}^{p}(\Omega_{t_{k}})$, $k=0,1,2,\dots,N-1$ are given.
The collection of these processes is denoted by $M_{\mathcal{E}}^{p,0}(0,T)$.

\begin{definition}
\label{d-Def-5 copy(1)}For an $\eta \in M_{\mathcal{E}}^{p,0}(0,T)$ with $\eta_{t}%
(\omega)=\sum_{k=0}^{N-1}\xi_{k}(\omega)\mathbf{1}_{[t_{k},t_{k+1})}(t)$, the
related integrals are
\begin{align*}
\int_{0}^{T}\eta_{t}(\omega)dt  &  :=\sum_{k=0}^{N-1}\xi_{k}(\omega
)(t_{k+1}-t_{k}),\\
\int_{0}^{T}\eta_{t}(\omega)dX_{t}  &  :=\sum_{k=0}^{N-1}\xi_{k}%
(\omega)(X_{t_{k+1}}-X_{t_{k}}),\\
\int_{0}^{T}\eta_{t}(\omega)dy_{t}  &  :=\sum_{k=0}^{N-1}\xi_{k}%
(\omega)(y_{t_{k+1}}-y_{t_{k}}).
\end{align*}

\end{definition}

\begin{definition}
$\label{mgp}$For each $p\geq1$, we denote by $M_{\mathcal{E}}^{p}(0,T)$ the completion
of $M_{\mathcal{E}}^{p,0}(0,T)$ under the norm%
\[
{{\left \Vert \eta \right \Vert _{M_{\mathcal{E}}^{p}(0,T)}:=\left \{  {\mathcal{E}%
}[{{\int_{0}^{T}|\eta_{t}|^{p}dt]}}\right \}  ^{1/p}}}.
\]

\end{definition}

It is clear that $M_{\mathcal{E}}^{p}(0,T)\supset M_{\mathcal{E}}^{q}(0,T)$ for $1\leq p\leq q.$
We also use $M_{\mathcal{E}}^{p}(0,T;\mathbb{R}^{d})$ for all $d$-dimensional stochastic
processes $\eta_{t}=(\eta_{t}^{1},\dots,\eta_{t}^{d})^T$, $t\geq0$ with
$\eta^{i}\in M_{\mathcal{E}}^{p}(0,T)$, $i=1,2,\dots,d$.

\subsubsection{$G$-It\^{o}'s calculus}

In the above space, it is easy to check that, for fixed $p_{i}\in \mathbb{R}^{2d}$, $i=1,\dots,n$,
the $n$-dimensional process defined by $(\left \langle Z_{t}(\omega
),p_{1}\right \rangle ,\cdots,\left \langle Z_{t}(\omega),p_{n}\right \rangle
)_{t\geq0}$ is also a Brownian motion, particularly $X_{t}+y_{t}$ is a
$G$-Brownian motion under $\mathcal{E}$.

Therefore, the process $z_{t}(\omega)=X_{t}(\omega
)+y_{t}(\omega)=\omega(t)$, $\omega \in \Omega=C([0,\infty),\mathbb{R}^{d})$
is a $G$-Brownian motion under $\mathcal{E}$, namely $(z_{t})_{t\geq0}$
is $c_{G}$-quasi surely continuous process such that the nonlinear distribution of
$z_{t}-z_{s}$ is that of $z_{t-s}$, and $z_{t}-z_{s}$ is independent from
$(z_{s_{1}},\dots,z_{s_{N}})$ for each $t\geq s\geq t_{s_{i}}$,
$i=1,2,\dots,N$.

It is worth noticing that $(z_{t})_{t\geq0}$ is also a nonlinear diffusion
process under $\widetilde{\mathcal{E}}$: for each $t_{i}\geq t\geq0$,
$i=1,\dots,N$, $\widetilde{\mathcal{E}}_{t}[\varphi(z_{t_{1}},\cdots,z_{t_{n}})]\text{ depends
only on }z_{t}$. Such nonlinear Markovian property plays an important role in this paper.

\begin{proposition}
If the function $\widetilde{G}$ is of the form $\widetilde{G}(x+y,p,A)$, then
$z_{t}=X_{t}+y_{t}$ still satisfies a martingale problem with nonlinear expectation
derived from the PDE
\[
\partial_{t}u(t,z)+\widetilde{G}(z,D_{z}u(t,z),D_{z}^2u(t,z))=0\quad (t,x)\in(0,\infty)\times \mathbb{R}^d.
\]
\end{proposition}

\noindent{\bf Proof.}
In this case the solution of the PDE%
\[
\partial_{t}u(t,x,y)-\widetilde{G}(x+y,D_{y}u(t,x,y),D_{x}^2%
u(t,x,y))=0,\  \  \ u(0,x,y)=\varphi(x+y)
\]
coincides with $\bar{u}(t,x+y)$, where $\bar{u}(t,z)$ is the solution to
the PDE
\begin{align*}
\partial_{t}\bar{u}(t,z)-G(z,D_{z}\bar{u}(t,z),D_{z}^2\bar{u}(t,z))  &
=0,\  \  \bar{u}(0,z)=\varphi(z),\  \\
(t,z)  &  \in(0,\infty)\times \mathbb{R}^{d}.
\end{align*}
\endofproof

Notice that nonlinear expectation $\widetilde{\mathcal{E}}$ is dominated by the sublinear expectation $\mathcal{E}$, the nonlinear expectation $\widetilde{\mathcal{E}}$ can still be defined on the Banach space $L^p_{{\mathcal{E}}},p\ge 1$. We give It\^{o}'s formula for a ``$G$-It\^{o}
process". For simplicity, we first consider the case of the function $\Phi$ being sufficiently regular and consider the general $n$-dimensional $G$-It\^{o}'s process%

\begin{equation}
\xi_{t}=\xi_{0}+\int_{0}^{t}\alpha_{s}ds+\int_{0}^{t}\beta_{s}dX_{s}+\int
_{0}^{t}\eta_{s}d\left \langle X\right \rangle _{s}+\int_{0}^{t}\kappa_{s}dy_{s},
\label{ItoProc}%
\end{equation}
where $\xi_0\in\mathbb{R}^n, \alpha_s \in \mathbb{R}^n, \beta_s,\kappa_s\in L(\mathbb{R}^d;\mathbb{R}^n),$ and $\eta_s\in L(\mathbb{R}^{d\times d};\mathbb{R}^n)$.

\begin{theorem}[It\^{o}'s formula]
\label{Thm6.5} Let $\Phi$ be a $C^{2}$-function on $\mathbb{R}^{n}$ such that
$\partial_{x^{\mu}x^{\nu}}^{2}\Phi$ satisfies the polynomial growth condition for
$\mu,\nu=1,\cdots,n$. Let $\alpha^{\nu}$, $\beta^{\nu j}$, and $\eta^{\nu}_{ ij}$,
$\nu=1,\dots,n$, $i,j=1,\dots,d$ be bounded processes in $M_{\mathcal{E}}^{2}(0,T)$.
Then for each $t\geq s \ge 0$ we have in $L_{\mathcal{E}}^{2}(\Omega_{t})$%
\begin{equation*}
\begin{aligned}
  \Phi(\xi_{t})-\Phi(\xi_{s})&=\int_{s}^{t}\partial_{x^{\nu}}\Phi(\xi_{\theta}%
)d\xi_{\theta}^\nu+\frac{1}{2}\sum_{i,j}\int_{s}^{t}\partial_{x^{\mu}x^{\nu}}^{2}\Phi(\xi
_{\theta})\beta_{\theta}^{\mu i}\beta_{\theta}^{\nu j}d\left \langle X \right \rangle^{ij}
_{\theta}\\
&= \int_s^t \langle D_x\Phi(\xi_\theta), d\xi_\theta\rangle + \frac{1}{2}\int_s^t\text{\em tr}[\beta^T D_x^2\Phi(\xi_\theta)\beta\,d\left \langle X\right \rangle_\theta],
\end{aligned}
\end{equation*}
where $\left \langle X \right \rangle^{ij} =\left \langle
X^{i},X^{j}\right \rangle.$
\end{theorem}

We have the following lemma.

\begin{proposition}
\label{Quadra-xi}
\[
\left \langle \xi \right \rangle _{t}\triangleq \left( \left \langle \xi \right \rangle_t^{ij}\right) = \left(\int_0^t \sum_{\mu,\nu}\beta^{i\mu}_s \beta^{j\nu}_s\,d\left\langle X^\mu, X^\nu\right\rangle _s\right) = \int_{0}^{t}\beta_{s}\,d\left \langle X\right \rangle _{s}\beta_{s}^{T}.%
\]

\end{proposition}

\noindent{\bf Proof.}
Let $\pi_{t}^{N}$, $N\in\mathbb{N}$, be a sequence of partitions of $[0,t]$.
We denote $\xi_{t}^{(N)}=\sum_{j=0}^{N-1}\xi_{t_{j}^{N}}1_{[t_{j},t_{j+1}%
)}(t)$, then
\begin{align*}
\xi_{t}\xi_t^{T}-\xi_{0}\xi_0^{T}  &  =\sum_{j=0}^{N-1}(\xi_{t_{j+1}^{N}}\xi_{t_{j+1}^{N}}^{T}-\xi
_{t_{j}^{N}}\xi
_{t_{j}^{N}}^{T})\\
&  =\int_{0}^{t}\left[\xi_{t}^{(N)}d\xi_{t}^T + d\xi_{t}\left(\xi_{t}^{(N)}\right)^T\right]+\sum_{j=0}^{N-1}(\xi_{t_{j+1}^{N}}%
-\xi_{t_{j}^{N}})(\xi_{t_{j+1}^{N}}%
-\xi_{t_{j}^{N}})^T.
\end{align*}
As $\mu(\pi_{t}^{N})\rightarrow0$, the first term on the right hand side
converges to ${{{\int_{0}^{t}\left[\xi_{s}\,d\xi^T_{s}+ d\xi_s\,\xi_s^T\right]}}}${{{ in }}}$L_{\mathcal{E}}%
^{2}(\Omega_{t})${, t{{he second one must be convergent, and we denote its limit
by $\left \langle \xi \right \rangle _{t}$, }}}%
\[
{{{\left \langle \xi \right \rangle _{t}=}}}\xi_{t}\xi_t^{T}-\xi_{0}\xi_0^{T} - \int_{0}^{t}\left[\xi_{s}\,d\xi^T_{s}+ d\xi_s\,\xi_s^T\right]%
\]
But by Theorem~\ref{Thm6.5}%
\[
\xi_{t}^{i}\xi_{t}^{j} - \xi_{0}^{i}\xi_{0}^{j} = \int_{0}^{t}(\xi^i_{\theta}\,d\xi_{\theta}^j +d\xi^i_{\theta}\,\xi_{\theta}^j)+ \sum_{\mu,\nu}\int_{0}^{t}\beta^{i\mu}_{s}\beta_{s}^{j\nu}d\left \langle
X\right \rangle ^{\mu\nu}_{s},\quad 1\le i,j\le n.
\]
\endofproof

\subsection{Proof of Theorems \ref{comparison3} and \ref{domthm}}



\paragraph{Proof of Theorem \ref{comparison3}.} The proof is based on  a key lemma.

\begin{lemma}\label{thm2.3}
  Suppose that each continuous function $G_i:[0,\infty)\times{\mathbb R}^d\times{\mathbb R}\times{\mathbb R}^d\times\mathbb{S}^d\to {\mathbb R}, i=0,1,\dots,k$ satisfies
  \begin{eqnarray}
  \lambda G_i(t,x,v,p,A) = G_i^{(1)}(t,x,v,\lambda p,\lambda A) + \lambda G_i^{(2)}(t,x,v,p),\label{homo}
  \end{eqnarray}
  for $\forall (t,x,v,p,A)\in [0,\infty)\times{\mathbb R}^d\times{\mathbb R}\times{\mathbb R}^d\times\mathbb{S}^d, \lambda \ge 0, i = 0,1,\dots,k$, where its respective decomposition functions $G_i^{(1)}$ and $G_i^{(2)}$ being continuous and satisfying the following condition
\begin{enumerate}[(A).]
  \setcounter{enumi}{6}
  \item (Condition (G))

  $|G(t,x,v,p,A) - G(t,y,v,p,A)| \le \bar{\rho}_G(1+(T-t)^{-1}+|x|+|y|+|v|)\rho_G(|x-y|+|p|\cdot|x-y|)$,\\for each $t\in[0,\infty), v\in{\mathbb R},x,y,p\in{\mathbb R}^d$, and $A\in\mathbb{S}^d$, where $\rho_G,\bar{\rho}_G:[0,\infty)\to[0,\infty)$ are continuous functions that satisfy $\rho_G(0)=0,\bar{\rho}_G(0)=0$.\label{g}
\end{enumerate}
Moreover, assume that
  \begin{equation*}
  \sum_{i=1}^k G_i(t,x,v_i,p_i,A_i) \le G_0(t,x,\sum_{i=1}^k v_i,\sum_{i=1}^k p_i,\sum_{i=1}^k A_i),
  \end{equation*}
  for each $(t,x,v_i,p_i,A_i) \in [0,\infty)\times{\mathbb R}^d\times{\mathbb R}\times{\mathbb R}^d\times\mathbb{S}^d$, and that
  \begin{align*}
    &G_0(t,x,v,p,A) \le G_0(t,x,v,p,A+\bar{A}),\\
    &|G_0(t,x,u,p,A) - G_0(t,x,v,p,\bar{\bar{A}})| \le L_0 (|u-v|+|A-\bar{\bar{A}}|),
  \end{align*}
  where $(t,x,p) \in [0,\infty)\times{\mathbb R}^d\times{\mathbb R}^d, u,v \in {\mathbb R}, A,\bar{A},\bar{\bar{A}}\in \mathbb{S}^d$ with $\bar{A}\ge0$, and $L_0>0$ is a constant. For each $i=1,\dots,k$, let $u_i \in \text{USC}([0,T]\times {\mathbb R}^d)$ be a viscosity subsolution of PDE
  \begin{equation*}
  \partial_t u(t,x) - G_i(t,x,u(t,x),Du(t,x),D^2u(t,x)) = 0, (t,x)\in(0,T]\times{\mathbb R}^d,
  \end{equation*}
  and let $u_0\in\text{LSC}([0,T]\times{\mathbb R}^d)$ be a viscosity supersolution of PDE
 \begin{equation*}
  \partial_t u(t,x) - G_0(t,x,u(t,x),Du(t,x),D^2u(t,x)) = 0, (t,x)\in(0,T]\times{\mathbb R}^d,
  \end{equation*}
  such that each $u_i, i = 0,1,\dots,k$ is with polynomial growth. Then
  $\displaystyle\sum_{i=1}^k u_i \le u_0$ on $(0,T]\times{\mathbb R}^d$ provided that $\displaystyle\sum_{i=1}^k u_i(0,x) \le u_0(0,x), x \in {\mathbb R}^d$.
\end{lemma}
Note that the above lemma  corrects  Theorem C.2.3 of \cite{P10a} where a condition of  type
(\ref{homo}) was missing.
Now, take $G_0 = G$, $G_1 = {\widetilde{G}}$, and define $G_2(t,x,p,A) = - {\widetilde{G}}(t,x,-p,-A)$.  Since ${\widetilde{G}}$ satisfies condition~\eqref{DOM}, we have
\begin{align*}
  G_1(t,x,p_1,A_1) + G_2(t,x,p_2,A_2) &= {\widetilde{G}}(t,x,p_1,A_1) - {\widetilde{G}}(t,x,-p_2,-A_2)\\
  &\le G(t,x,p_1- (-p_2),A_1 - (-A_2))\\
  & = G_0(t,x,p_1+p_2, A_1+A_2).
\end{align*}
Applying Lemma~\ref{thm2.3} yields Theorem \ref{comparison3}.
\endofproof

\paragraph{Proof of Theorem \ref{domthm}}
    Let $\Phi\in C^2([0,T]\times{\mathbb R}^d)$, and $w-\Phi$ achieves its global maximum at $(t_0,x_0)\in (0,T)\times{\mathbb R}^d$. Set
    \[
    \Psi_{{\varepsilon},\delta}(t,x,s,y) = u_1(t,x) - u_2(s,y) - \frac{|x-y|^2}{{\varepsilon}^2} - \frac{|t-s|^2}{\delta^2} - \Phi(t,x),
    \]
    where ${\varepsilon}, \delta>0$.
    Since $(t_0,x_0)$ is a strict global maximum point of $w-\Phi$, as the proof for Lemma 3.1 in \cite{CIL92}, there exists a sequence $(\tilde{t},\tilde{x},\tilde{s},\tilde{y}) = (\tilde{t}({\varepsilon},\delta),\tilde{x}({\varepsilon},\delta),\tilde{s}({\varepsilon},\delta),\tilde{y}({\varepsilon},\delta))$ such that
    \begin{itemize}
      \item $(\tilde{t},\tilde{x},\tilde{s},\tilde{y})$ is a global maximum point of $\Psi_{{\varepsilon},\delta}$ in $([0,T]\times\bar{B}_R)^2$;
      \item $(\tilde{t},\tilde{x}),(\tilde{s},\tilde{y})\to (t_0,x_0)$  as $({\varepsilon},\delta)\to 0$;
      \item $\dfrac{|\tilde{x}-\tilde{y}|}{{\varepsilon}^2}$ and $\dfrac{|\tilde{t}-\tilde{s}|}{\delta^2}$ are bounded and tend to zero as $({\varepsilon},\delta)\to 0$.
    \end{itemize}
    According to  Theorem 8.3 in~\cite{CIL92}, there exist $X, Y \in \mathbb{S}^d$ such that
    \begin{align*}
      \left(\frac{2(\tilde{t}-\tilde{s})}{\delta^2}+\frac{\partial \Phi}{\partial t}(\tilde{t},\tilde{x}), \frac{2(\tilde{x}-\tilde{y})}{{\varepsilon}^2}+D\Phi(\tilde{t},\tilde{x}),X\right) &\in \bar{P}^{2,+}_{B_R}u_1(\tilde{t},\tilde{x}),\\
      \left(\frac{2(\tilde{t}-\tilde{s})}{\delta^2}, \frac{2(\tilde{x}-\tilde{y})}{{\varepsilon}^2},Y\right) &\in \bar{P}^{2,-}_{B_R}u_2(\tilde{s},\tilde{y}),\\
      \left(\begin{matrix} X & 0\\ 0 & -Y\end{matrix}\right) \le \frac{4}{{\varepsilon}^2}\left(\begin{matrix}I&-I\\-I&I\end{matrix}\right)&+ \left(\begin{matrix}D^2\Phi(\tilde{t},\tilde{x})&0\\0&0\end{matrix}\right).
    \end{align*}

    Without loss of generality, assume that $(\tilde{t},\tilde{x},\tilde{s},\tilde{y})$ is a global maximum point of $\Psi_{{\varepsilon},\delta}$ in $([0,T]\times{\mathbb R}^d)^2$. Since $u_1$ and $u_2$ are subsolution and supersolution to the PDE ~\eqref{GE2}, we have
   \begin{equation*}
      -\frac{2(\tilde{t}-\tilde{s})}{\delta^2} - \frac{\partial \Phi}{\partial t}(\tilde{t},\tilde{x}) - {\widetilde{G}}(\tilde{x},\frac{2(\tilde{x}-\tilde{y})}{{\varepsilon}^2}+D\Phi(\tilde{t},\tilde{x}),X) \le 0
    \end{equation*}
    and
    \begin{equation*}
      -\frac{2(\tilde{t}-\tilde{s})}{\delta^2}  - {\widetilde{G}}(\tilde{y},\frac{2(\tilde{x}-\tilde{y})}{{\varepsilon}^2}),Y) \ge 0.
    \end{equation*}
    Therefore
    \begin{eqnarray}\label{key1}
    - \frac{\partial \Phi}{\partial t}(\tilde{t},\tilde{x}) - \left[{\widetilde{G}}(\tilde{x},\frac{2(\tilde{x}-\tilde{y})}{{\varepsilon}^2}+D\Phi(\tilde{t},\tilde{x}),X)-{\widetilde{G}}(\tilde{y},\frac{2(\tilde{x}-\tilde{y})}{{\varepsilon}^2},Y)\right]\le 0.
    \end{eqnarray}
    Note that
    \begin{align*}
      &{\widetilde{G}}(\tilde{x},\frac{2(\tilde{x}-\tilde{y})}{{\varepsilon}^2}+D\Phi(\tilde{t},\tilde{x}),X)-{\widetilde{G}}(\tilde{y},\frac{2(\tilde{x}-\tilde{y})}{{\varepsilon}^2},Y)\\
      =\,&\sup_{{\Gamma}}\inf_{{\Lambda}}\left\{{\frac{1}{2}}\text{tr}[{\sigma}(\tilde{x},{\gamma},{\lambda}){\sigma}^T(\tilde{x},{\gamma},{\lambda})X] + \left\langle b(\tilde{x},{\gamma},{\lambda}),\frac{2(\tilde{x}-\tilde{y})}{{\varepsilon}^2}+D\Phi(\tilde{t},\tilde{x})\right\rangle\right\}\\
      \, &- \sup_{{\Gamma}}\inf_{{\Lambda}}\left\{{\frac{1}{2}}\text{tr}[{\sigma}(\tilde{y},{\gamma},{\lambda}){\sigma}^T(\tilde{y},{\gamma},{\lambda})Y] + \left\langle b(\tilde{x},{\gamma},{\lambda}),\frac{2(\tilde{x}-\tilde{y})}{{\varepsilon}^2}\right\rangle\right\}\\
      \le\,&\sup_{{\Gamma},{\Lambda}}\bigg\{{\frac{1}{2}}\text{tr}[{\sigma}(\tilde{x},{\gamma},{\lambda}){\sigma}^T(\tilde{x},{\gamma},{\lambda})X-{\sigma}(\tilde{y},{\gamma},{\lambda}){\sigma}^T(\tilde{y},{\gamma},{\lambda})Y]+\left\langle b(\tilde{x},{\gamma},{\lambda})-b(\tilde{y},{\gamma},{\lambda}),\frac{2(\tilde{x}-\tilde{y})}{{\varepsilon}^2}\right\rangle\\ \,&\quad\quad+ \langle D\Phi(\tilde{t},\tilde{x}),b(\tilde{x},{\gamma},{\lambda})\rangle\bigg\},
      \end{align*}
      \begin{align*}
      &\text{tr}[{\sigma}(\tilde{x},{\gamma},{\lambda}){\sigma}^T(\tilde{x},{\gamma},{\lambda})X-{\sigma}(\tilde{y},{\gamma},{\lambda}){\sigma}^T(\tilde{y},{\gamma},{\lambda})Y]\\
            =\,&\text{tr}\left[\left(\begin{matrix}{\sigma}(\tilde{x},{\gamma},{\lambda}){\sigma}^T(\tilde{x},{\gamma},{\lambda})&{\sigma}(\tilde{x},{\gamma},{\lambda}){\sigma}^T(\tilde{y},{\gamma},{\lambda})\\{\sigma}(\tilde{y},{\gamma},{\lambda}){\sigma}^T(\tilde{x},{\gamma},{\lambda})&{\sigma}(\tilde{y},{\gamma},{\lambda}){\sigma}^T(\tilde{y},{\gamma},{\lambda})\end{matrix}\right)\left(\begin{matrix}X&0\\0&-Y\end{matrix}\right)\right]\\
      \le\,&\frac{4}{{\varepsilon}^2}\text{tr}\left[\left(\begin{matrix}{\sigma}(\tilde{x},{\gamma},{\lambda}){\sigma}^T(\tilde{x},{\gamma},{\lambda})&{\sigma}(\tilde{x},{\gamma},{\lambda}){\sigma}^T(\tilde{y},{\gamma},{\lambda})\\{\sigma}(\tilde{y},{\gamma},{\lambda}){\sigma}^T(\tilde{x},{\gamma},{\lambda})&{\sigma}(\tilde{y},{\gamma},{\lambda}){\sigma}^T(\tilde{y},{\gamma},{\lambda})\end{matrix}\right)\left(\begin{matrix}I&-I\\-I&I\end{matrix}\right)\right]\\ \,&+ \text{tr}[{\sigma}(\tilde{x},{\gamma},{\lambda}){\sigma}^T(\tilde{x},{\gamma},{\lambda})D^2\Phi(\tilde{t},\tilde{x})]\\
      =\,&\frac{4}{{\varepsilon}^2}\text{tr}[({\sigma}(\tilde{x},{\gamma},{\lambda}){\sigma}^T(\tilde{x},{\gamma},{\lambda})-{\sigma}(\tilde{x},{\gamma},{\lambda}){\sigma}^T(\tilde{y},{\gamma},{\lambda}))({\sigma}(\tilde{y},{\gamma},{\lambda}){\sigma}^T(\tilde{x},{\gamma},{\lambda})-{\sigma}(\tilde{y},{\gamma},{\lambda}){\sigma}^T(\tilde{y},{\gamma},{\lambda}))]\\
      \,&+ \text{tr}[{\sigma}(\tilde{x},{\gamma},{\lambda}){\sigma}^T(\tilde{x},{\gamma},{\lambda})D^2\Phi(\tilde{t},\tilde{x})]\\
      =\,&\frac{4}{{\varepsilon}^2}\text{tr}[({\sigma}(\tilde{x},{\gamma},{\lambda})-{\sigma}(\tilde{y},{\gamma},{\lambda}))({\sigma}^T(\tilde{x},{\gamma},{\lambda})-{\sigma}^T(\tilde{y},{\gamma},{\lambda}))]+ \text{tr}[{\sigma}(\tilde{x},{\gamma},{\lambda}){\sigma}^T(\tilde{x},{\gamma},{\lambda})D^2\Phi(\tilde{t},\tilde{x})]\\
      \le\,& 4\tilde{L}^2\frac{|\tilde{x}-\tilde{y}|^2}{{\varepsilon}^2} + \text{tr}[{\sigma}(\tilde{x},{\gamma},{\lambda}){\sigma}^T(\tilde{x},{\gamma},{\lambda})D^2\Phi(\tilde{t},\tilde{x})],
     \end{align*}
    and
    \[
    \left\langle b(\tilde{x},{\gamma},{\lambda})-b(\tilde{y},{\gamma},{\lambda}),\frac{2(\tilde{x}-\tilde{y})}{{\varepsilon}^2}\right\rangle \le\tilde{L}\frac{|\tilde{x}-\tilde{y}|^2}{{\varepsilon}^2}.
    \]
    We have
    \begin{eqnarray}\label{key2}
    - \frac{\partial \Phi}{\partial t}(\tilde{t},\tilde{x}) - G(\tilde{x},D\Phi(\tilde{t},\tilde{x}),D^2\Phi(\tilde{t},\tilde{x}))\le 5\tilde{L}\frac{|\tilde{x}-\tilde{y}|^2}{{\varepsilon}^2},
    \end{eqnarray}
    and the right hand side tends to 0 as $({\varepsilon},\delta)\to 0$.
    \endofproof




\begin{thebibliography}{99}

\bibitem[ADEH99]{ADEH99}
{\sc P.\,Artzner, F.\,Delbaen, J.-M.\,Eber, and D.\,Heath} (1999), Coherent measures of risk, {\em Mathematical Finance}, Vol. 9, 203-228.

\bibitem[BBP97]{BBP97}
{\sc G.\,Barles, R.\,Buckdahn, and E.\,Pardoux} (1997), Backward stochastic differential equations and integral-partial differential equations, {\em Stochasitcs and Stochastics Reports\/}, Vol. 60, 57-83.

\bibitem[BL08]{BL08}
{\sc R.\,Buckdahn and J.\,Li} (2008), Stochastic differential games and viscosity solutions of Hamilton-Jacobi-Bellman-Isaacs equations, {\em SIAM Journal on Control and Optimization\/}, Vol. 47, No. 1, 444-475.


\bibitem[CSTV07]{CSTV07}
{\sc P.\,Cheridito, H.\,M.\,Soner, N.\,Touzi, and N.\,Victoir} (2007), Second-order backward stochastic differential equations and fully
nonlinear parabolic PDEs, {\em Communications on Pure and Applied Mathematics}, Vol. LX, 1081-1110.



\bibitem[CIL92]{CIL92}
{\sc M.\,G.\,Crandall, H.\,Ishii, and P.\,L.\,Lions} (1992), User's guide to viscosity solutions of second order partial differential equations, {\em Bulletin of The American Mathematical Society\/}, Vol. 27, No.1, 1-67.

\bibitem[DHP11]{DHP11}
{\sc L.\,Denis, M.\,S.\,Hu, and S.\,G.\,Peng} (2011), Function spaces and capacity related to a sublinear expectation: application to $G$-Brownian motion paths, {\em Potential Analysis\/}, Vol. 34, Issue 2, 139--161.


\bibitem[EK05]{EK05}
{\sc S.\,N.\,Ethier and T.\,G.\,Kurtz} (2005), {\em Markov Processes: characterization and convergence}, J. Wiley and Sons, New York.


\bibitem[FS92]{FS92}
{\sc W.\,H.\,Fleming and H.\,M.\,Soner} (1992), {\em Controlled Markov Processes and Viscosity Solutions}, Springer-Verlag, New York.

\bibitem[HP09]{HP09}
{\sc M.\,Hu and S.\,G.\,Peng} (2009), $G$-L\'{e}vy processes under sublinear expectations,  {\em Preprint, arXiv:0911.3533v1 [math.PR] 18 Nov 2009.\/}

\bibitem[KS91]{KS91}
{\sc I.\,Karatzas and S.\,E.\,Shreve} (1991), {\em Brownian Motion and Stochastic Calculus\/}, 2nd edition, Springer-Verlag, New York.


\bibitem[K87]{K87}
{\sc N.\,V.\,Krylov} (1987),
{\em Nonlinear Elliptic and Parabolic Equations of the Second Order\/}, D.\,Reidel Publishing Company (Original Russian version by Nauka, Moscow, 1985).

\bibitem[LP11]{LP11}
{\sc X.\,P.\,Li and S.\,G.\,Peng} (2011), Stopping times and related It\^{o}'s calculus with $G$-Brownian motion, {\em Stochastic Processes and Applications\/}, Vol. 121, Issue 7, 1492-1508.

\bibitem[N76a]{N76a}
 {\sc M.\,Nisio} (1976), On a nonlinear semigroup attached to optimal stochastic control. {\em Publications of the Research Institute for Mathematical Sciences}, 12(2): 513-537.

\bibitem[N76b]{N76b}
 {\sc M.\,Nisio} (1976), On stochastic optimal controls and envelope of Markovian semigroups. {\em Proceedings of International Symposium, Kyoto}, 297-325.

\bibitem[N13]{N12}
{\sc M.\,Nutz} (2013), Random $G$-expectations,  {\em The Annals of Applied Probability\/}, Vol. 23, No. 5, 1755-1777.



\bibitem[{\O}03]{O03}
{\sc B.\,{\O}ksendal} (2003),
{\em Stochastic Differential Equations. An Introduction with
Applications\/}, 6th edition, Springer-Verlag, New York.



  \bibitem[O11]{O11}
  {\sc E.\,Osuka} (2011), Girsanov's theorem for $G$-Brownian motion, {\it Preprint, arXiv:1106.2387v1.}

\bibitem[P05]{P05}
{\sc S.\,G.\,Peng} (2005), Nonlinear expectations and nonlinear Markov chains, {\em Chinese Annals of Mathematics\/}, Vol. 26, No.2, 159-184.

\bibitem[P07a]{P07a}
{\sc S.\,G.\,Peng} (2007), $G$-expectation, $G$-Brownian motion and related stochastic calculus of It\^o type, {\em Stochastic Analysis and Applications, The Abel Symposium 2005, Abel Symposia, Vol.2,\/}  541-567, Springer-Verlag, Berlin Heidelberg.

\bibitem[P07b]{P07b}
{\sc S.\,G.\,Peng} (2007), {\em Lecture Notes: $G$-Brownian Motion and Dynamic Risk Measure under Volatility Uncertainty}, Preprint, arXiv:0711.2834v1.


\bibitem[P10a]{P10a}
{\sc S.\,G.\,Peng} (2010), {\em Nonlinear Expectations and Stochastic Calculus under Uncertainty}, Preprint, arXiv:1002.4546v1.

\bibitem[P10b]{P10b}
{\sc S.\,G.\,Peng} (2010), Tightness, weak compactness of nonlinear expectations and application to CLT, {\it Preprint, arXiv:1006.2541v1.}

\bibitem[Ph09]{Ph09}
{\sc H.\,Pham} (2009), {\em Continuous-time Stochastic Control and Optimization with Financial Applications}, Springer-Verlag, New York.

\bibitem[PSZ12]{PSZ12}
{\sc S.\,G.\,Peng, Y.\,S.\,Song, and J.\,F.\,Zhang} (2012), A complete representation theorem for $G$-martingales, {\it Preprint, arXiv:1201.2629v1.}

\bibitem[STZ11a]{STZ11a}
{\sc H.\,M.\,Soner, N.\,Touzi, and J.\,F.\,Zhang} (2011), Martingale representation theorem for the $G$-expecation, {\em Stochastic Processes and their Applications\/}, Vol. 121, Issue 2, 265-287.

\bibitem[STZ11b]{STZ11b}
{\sc H.\,M.\,Soner, N.\,Touzi, and J.\,F.\,Zhang} (2011), Quasi-sure stochastic analysis through aggregation, {\em Electronic Journal of Probability}, Vol. 16, No. 67, 1844-1879.

\bibitem[STZ12]{STZ12}
{\sc H.\,M.\,Soner, N.\,Touzi, and J.\,F.\,Zhang} (2012), Well-posedness of second order backward SDEs,  {\em Probability Theory and Related Fields\/}, Vol. 153, Issue 1-2, 149-190.


\bibitem[SV69]{SV69}
 {\sc D.W.\,Stroock and S.R.S.\,Varadhan} (1969), Diffusion processes with continuous coefficients, I and II, {\em Communications on Pure and Applied Mathematics},
    Vol. 22, 345-400 and 479-530.

\bibitem[SV79]{SV79}
{\sc D.W.\,Stroock and S.R.S.\,Varadhan} (1979), {\em Multidimensional Diffusion Processes\/}, Springer, New York.


\bibitem[XZ10]{XZ10}
{\sc J.\,Xu and B.\,Zhang} (2010), Martingale property and capacity under $G$-framework, {\em Electronic Journal of Probability\/}, Vol. 15, No. 67, 2041-2068.

\bibitem[XHZ11]{XHZ11}
{\sc J.\,Xu, H.\,Shang, and B.\,Zhang} (2011), A Girsanov type theorem under $G$-framework, {\em Stochastic Analysis and Applications\/}, Vol. 29, Issue 3, 386-406.

\bibitem[YZ99]{YZ99}
{\sc J.\,M.\,Yong and X.\,Y.\,Zhou} (1999), {\em Stochastic Controls: Hamiltonian Systems and HJB Equations\/}, Springer-Verlag, New York.

\end{thebibliography}

\end{document}